\documentclass[11pt,reqno]{amsart}
\usepackage[utf8]{inputenc}
\usepackage{amsmath}
\usepackage{amsfonts}
\usepackage{amssymb}
\usepackage{amsthm}
\usepackage{appendix}
\usepackage{indentfirst}
\usepackage{mathtools}
\mathtoolsset{showonlyrefs=true}
\usepackage{enumitem}
\usepackage{hyperref}

\hypersetup{
colorlinks=true,
linkcolor=blue,
citecolor=blue,
}

\title[Profiled expansive solutions]{Rigidity and existence of profiled expansive solutions of the $N$-body problem with homogeneous potentials}
\author{Yutong Wu}
\address{Department of Mathematics, Yale University, New Haven, CT 06511} 
\email{yutong.wu.yw894@yale.edu}

	

\def \N{\mathbb{N}}
\def \R{\mathbb{R}}

\def \d{\mathrm{d}}
\def \dt{\mathrm{d}t}
\def \X{\mathcal{X}}
\def \Y{\mathcal{Y}}
\def \D{\mathcal{D}}
\def \A{\mathcal{A}}
\def \L{\mathcal{L}}

\def \H{\mathcal{H}}

\newcommand{\f}[2]{\frac{#1}{#2}}

\begin{document}

\numberwithin{equation}{section} 

\newtheorem{thm}{Theorem}
\newtheorem{lem}{Lemma}[section]
\newtheorem{cor}[thm]{Corollary}
\newtheorem{prop}[lem]{Proposition}

\theoremstyle{remark}
\newtheorem{rmk}[lem]{Remark}

\theoremstyle{definition}
\newtheorem{defn}[lem]{Definition}

\begin{abstract}
We study expansive solutions of the $N$-body problem with homogeneous potentials of degree $-p$, where $p>0$. We establish asymptotic rigidity for arbitrary solutions with admissible profiles. For general admissible profiles and $0<p<2$, we obtain convergent power-log expansions whose exponents are governed by the leading profiles. In the hyperbolic case, we derive a more explicit analytic representation for all $p>0$, extending Chazy's Newtonian result. Together with recent results in Wu~\cite{Wu}, this yields a bijective parametrization of hyperbolic solutions near infinity by their asymptotic velocities and offsets.

We also establish the existence of free-time minimizing solutions with arbitrary prescribed initial configurations and admissible asymptotic profiles for $0<p<2$, covering the hyperbolic, parabolic and hyperbolic-parabolic regimes. In the hyperbolic case, the result extends to all $p>0$, including possibly collisional initial configurations in the strong-force regime $p \ge 2$. Furthermore, we prove that when $p \neq 2$, every solution with an admissible asymptotic profile is eventually a free-time minimizer.
\end{abstract}

\maketitle

\section*{Declaration of AI use}
Generative AI tools were used for literature search and proofreading. In particular, in Lemma~\ref{lem generators}, the author originally proved the vanishing of resonant terms for $p \in \N_{>1}$ by comparing coefficients, using an argument analogous to that in Lemma~\ref{lem low order hyperbolic}. AI subsequently suggested an alternative approach using the normal form structure directly. This argument was then completed and adopted by the author. Apart from this, the mathematical ideas, results, and proofs presented in the paper were independently developed by the author.

\section{Introduction}

We study the $N$-body problem with a homogeneous potential. Given the dimension $d \in \N_{\ge 1}$, the number of bodies $N \in \N_{\ge 2}$, the masses of the bodies $\lambda_1, \cdots, \lambda_N>0$ and the homogeneous exponent $p>0$, we consider the following system
\begin{equation} \label{eq N-body}
    \lambda_i \ddot{x}_i= -\sum_{j \neq i} \frac{\lambda_i \lambda_j(x_i- x_j)}{|x_i- x_j|^{p+2}}, \quad \forall i=1,2, \cdots, N,
\end{equation}
where the position $x_i$ takes values in $\R^d$ for each $1 \le i \le N$. 

We can write \eqref{eq N-body} in a more compact form. Set
\begin{equation}
    \X= \R^{dN}, \quad \Y= \left\{ x \in \X \ \big| \ x_j \neq x_k, \ \forall j \neq k \right\}, \quad \Delta= \X \setminus \Y.
\end{equation}
For $x= (x_1, \cdots, x_N) \in \X$, define
\begin{equation}
    U(x):= \frac{1}{p} \sum_{i<j} \frac{\lambda_i \lambda_j}{|x_i- x_j|^p} \quad \text{and} \quad |x|^2 _{\lambda}:= \sum_{i=1}^N \lambda_i |x_i|^2.
\end{equation}
We see $U(x)= +\infty \iff x \in \Delta$. Besides, $| \cdot |^2_{\lambda}$ induces a metric on $\X$. Let $\nabla$ denote the gradient with respect to this metric. In other words, $\nabla_i= \frac{1}{\lambda_i} \partial_{x_i}$. Then \eqref{eq N-body} is equivalent to 
\begin{equation} 
    \ddot{x}= \nabla U(x).
\end{equation}

We focus on the \emph{expansive solutions} of \eqref{eq N-body}, namely those for which all mutual distances diverge to infinity:
\begin{equation}
    |x_i(t)- x_j(t)| \to +\infty \quad \text{as} \;\; t \to +\infty, \quad \forall i\neq j.
\end{equation}

The asymptotic structure of expansive solutions was studied in Marchal--Saari \cite{MarchalSaari}, who proved the existence of limiting velocities, and further established parabolic growth of internal motions in the Newtonian case $p=1$. Recently, Wang--Yu \cite{WangYu} proved the existence of a limiting central configuration within each cluster when $p=1$. Motivated by these results, we consider expansive solutions with such limiting profiles for general homogeneous potentials. 

\begin{defn}[Profiled expansive solution] \label{def profiled expansive}
We call a solution $x$ of \eqref{eq N-body} a \textbf{profiled expansive solution} if it admits the asymptotic expansion
\begin{equation}
    x(t)= vt+ at^{\frac{2}{p+2}}+ o(t^{\frac{2}{p+2}}) \quad \text{as} \;\; t \to +\infty
\end{equation}
for some $v, a \in \X$, where $v_i= v_j$ for $i \neq j$ implies $a_i \neq a_j$.
\end{defn}

Clearly, every profiled expansive solution is expansive. In the Newtonian case $p=1$, the converse also holds by \cite{MarchalSaari, WangYu}.

For $v=(v_1, \cdots, v_N) \in \X$, the relation $i \sim_v j$ given by $v_i= v_j$ is an equivalence relation on $\{1, \cdots, N\}$. Let $M_v$ denote the set of equivalence classes and we call each equivalence class a \emph{cluster}. For a profiled expansive solution, we have
\begin{equation}
    |x_i(t)- x_j(t)| \sim \begin{cases}
        t, & i \not\sim_v j, \\
        t^{\frac{2}{p+2}}, & i \sim_v j, \; i \neq j,
    \end{cases} \quad \text{as} \;\; t \to +\infty.
\end{equation}
Following Chazy's terminology, we distinguish three types of profiled expansive solutions.
\begin{itemize}
    \item \textbf{Hyperbolic:} $|M_v|=N$, which means each cluster is a singleton.
    \item \textbf{Parabolic:} $|M_v|= 1$, which means there is only one cluster.
    \item \textbf{Hyperbolic-parabolic:} $1< |M_v|< N$, which means there is more than one cluster, and at least one cluster contains more than one body.
\end{itemize}

Let $v \in \X$ and $I \in M_v$ be a cluster. For $a = (a_1, \cdots a_N) \in \X$, we let $a_I$ denote the projection of $a$ onto the subspace of $\X$ indexed by $I$. More precisely, we let $a_I= (a_{I,1}, \cdots a_{I,N})$, where 
\begin{equation}
    a_{I,i} = a_i \;\; \text{if} \;\; i \in I; \qquad a_{I,i} = 0 \;\; \text{if} \;\; i \not\in I.
\end{equation}
For each cluster $I \in M_v$, we define the internal potential
\begin{equation}
    U_I(x):= \frac{1}{p} \sum_{i<j \in I} \frac{\lambda_i \lambda_j}{|x_i- x_j|^p}.
\end{equation}
Clearly, $U_I(x)= U_I(x_I)$, and $U_I \equiv 0$ if $|I|=1$.

\begin{defn}[Admissible pair] \label{def admissible pair}
We call $(v, a) \in \X \times \X$ an \textbf{admissible pair} if for each cluster $I \in M_v$, the following conditions hold for $a_I$.
\begin{enumerate}
    \item \label{item minimality} \emph{Minimality:} $U_I(a_I)= \min \{ U_I(b_I) \ | \ b \in \X, |b_I|_\lambda= |a_I|_\lambda \}$;
    \item \emph{Normalization:} $\nabla U_I(a_I) = -\frac{2p}{(p+2)^2} a_I$.
\end{enumerate}
\end{defn}

For a cluster $I$ with $|I| \ge 2$, a configuration satisfying (1) can be rescaled to satisfy (2). The particular constant in (2) ensures that $t \mapsto a_I t^{2/{(p+2)}}$ is a homothetic parabolic solution of the isolated $|I|$-body problem. For a singleton cluster $I=\{i\}$, (2) implies $a_i=0$.

The admissibility condition will play a central role in both our asymptotic rigidity analysis and the variational construction of expansive solutions. 
\;\\

Now we state the main results of this paper.

\subsection{Asymptotic rigidity} \

The first part of our main results concerns asymptotic rigidity. In the Newtonian case $p=1$, every expansive solution admits an asymptotic profile as in Definition~\ref{def profiled expansive}, by the results of Marchal--Saari~\cite{MarchalSaari} and Wang--Yu~\cite{WangYu}. This naturally raises the question of what further information can be obtained about the asymptotic behavior of such solutions, particularly possible correction orders and the structure of subleading terms. The problem becomes more delicate in the hyperbolic-parabolic regime, where intra-cluster and inter-cluster interactions have different scales. We will investigate these issues for arbitrary solutions with admissible asymptotic profiles.

The following result establishes convergent asymptotic expansions in all hyperbolic, parabolic and hyperbolic-parabolic regimes.

\begin{thm}[Rigidity of profiled expansive solutions] \label{thm rigidity mixed} \

Assume $d \in \N_{\ge 1}$ and $p \in (0,2)$. Let $(v,a)$ be an admissible pair. Assume $x \in C^2 \big( (T, +\infty); \Y \big)$ solves \eqref{eq N-body} on $(T,+\infty)$ for some $T>0$ and
\begin{equation} \label{eq expansion profiled expansive}
    x(t) = vt + at^{\frac{2}{p+2}} + o(t^{\frac{2}{p+2}}) \quad \text{as} \;\; t \to +\infty.
\end{equation}
Then there exist $K \in \N \cup \{\infty\}$, strictly decreasing real numbers $\{\alpha_i\}_{i=1}^K$ and nonzero $\X$-valued polynomials $\{H_i\}_{i=1}^K$ such that
\begin{equation} \label{eq all-order expansion of solution}
    x(t) = vt + at^{\frac{2}{p+2}}+ \sum_{i=1}^K t^{\alpha_i} H_i(\log t) \quad \text{for} \;\; t \gg 1,
\end{equation}
where the series converges absolutely and uniformly for $t$ large enough. Moreover, the exponent $\alpha_i$ belongs to a discrete subset of $\R$ that depends only on $a$ and the cluster partition $M_v$. We also have 
\begin{equation} \label{eq mixed exponent bound}
    \alpha_i \le \max \{\frac{p}{p+2}, 1-p \}
\end{equation} 
and 
\begin{equation} \label{eq mixed no resonant}
    \begin{cases}
         \alpha_i \ge 1-p \\
         \alpha_i> \frac{2p-2}{p+2}
    \end{cases} \implies H_i \; \text{is a constant}.
\end{equation}

In the parabolic case, namely if $|M_v|=1$, then \eqref{eq mixed exponent bound} can be strengthened to
\begin{equation}
    \alpha_i \le \frac{p}{p+2}
\end{equation}
and \eqref{eq mixed no resonant} can be strengthened to
\begin{equation}
    \alpha_i> \frac{2p-2}{p+2} \implies H_i \; \text{is a constant}.
\end{equation}
\end{thm}

\begin{rmk} \label{rmk sharp estimate}
By \eqref{eq mixed exponent bound} and \eqref{eq mixed no resonant}, we deduce that any solution $x$ satisfying \eqref{eq expansion profiled expansive} automatically admits the stronger remainder estimate
\begin{equation} \label{eq expansive solutoin sharp remainder}
    x(t) = vt + at^{\frac{2}{p+2}} + O(t^{\max\{ \frac{p}{p+2}, 1-p\}}) \quad \text{as} \;\; t \to +\infty.
\end{equation}
In the parabolic case ($|M_v|=1$), \eqref{eq expansive solutoin sharp remainder} can be further strengthened to
\begin{equation}
    x(t) = vt + at^{\frac{2}{p+2}} + O(t^{\frac{p}{p+2}}) \quad \text{as} \;\; t \to +\infty.
\end{equation}
\end{rmk}

In the hyperbolic case, the result holds for all $p>0$, and we have more precise information of the expansion.

\begin{thm}[Rigidity of hyperbolic solutions] \label{thm rigidity hyperbolic} \

Assume $d \in \N_{\ge 1}$ and $p>0$. Let $v \in \Y$. Assume $x \in C^2\big( (T,+\infty); \Y \big)$ solves \eqref{eq N-body} on $(T,+\infty)$ for some $T>0$ and one has
\begin{equation} \label{eq expansion hyperbolic}
    x(t)= vt+ o(t) \quad \text{as} \;\; t \to +\infty.
\end{equation}
Then there exists an $\X$-valued function $\Phi( \cdot, \cdot)$, which is analytic in a neighborhood of $(0,0)$ and satisfies $\Phi(0,0)=0$, such that
\begin{equation} \label{eq hyperbolic analytic function}
    x(t)= \begin{cases}
        vt+ t\Phi(t^{-p}, t^{-1}), & p^{-1} \not\in \N, \\
        vt+ t\Phi(t^{-p}, t^{-1} \log t), & p^{-1} \in \N
    \end{cases}
    \quad \text{for} \;\; t \gg 1.
\end{equation}
\end{thm}

\begin{rmk}
The analytic function representation \eqref{eq hyperbolic analytic function} yields an absolutely and uniformly convergent series expansion as in Theorem~\ref{thm rigidity mixed}
\begin{equation} 
    x(t) = vt+ \sum_{i=1}^K t^{\alpha_i} H_i(\log t) \quad \text{for} \;\; t \gg 1.
\end{equation}
Moreover, we have $\alpha_i \in \{ 1-m_1-m_2p \ | \ m_1, m_2 \in \N,\ m_1+m_2>0 \}$, $\alpha_i \le \max \{0, 1-p\}$ and $H_i$ is a constant whenever $\alpha_i>0$.

The expansion in Theorem~\ref{thm rigidity mixed} also admits an analytic expression in suitable power-log variables similar to \eqref{eq hyperbolic analytic function}, as follows from the proof. This is also how we obtain a discrete subset of $\R$ containing all possible exponents. However, these variables do not admit the same universal explicit description as in the hyperbolic case.
\end{rmk}

Combining Theorem~\ref{thm rigidity hyperbolic} with Wu's existence and uniqueness result for hyperbolic solutions with prescribed asymptotic data ~\cite[Theorem~1]{Wu}, we obtain a complete parametrization of hyperbolic solutions. 

\begin{cor}[Bijective parametrization of hyperbolic solutions] \label{cor hyperbolic} \

Assume $d \in \N_{\ge 1}$ and $p>0$. Let $n$ be the largest nonnegative integer such that $np<1$. Then there exist analytic functions $c_1, \cdots, c_{n+1}: \Y \to \X$ with the following properties. 

For any hyperbolic solution $x$ of \eqref{eq N-body}, there exist $v \in \Y$ and $b \in \X$ such that: if $p^{-1} \not\in \N$, then
\begin{equation} 
    x(t)= vt+ \sum_{j=1}^n c_j(v) t^{1-jp}+ b+ c_{n+1}(v) t^{1-(n+1)p}+ O(t^{-p}) \;\; \text{as} \;\; t \to +\infty.
\end{equation}  
If $p^{-1} \in \N$, then
\begin{equation}
    x(t)= vt+ \sum_{j=1}^n c_j(v) t^{1-jp}+ c_{n+1}(v) \log t+ b+ O(t^{-p} \log t) \;\; \text{as} \;\; t \to +\infty.
\end{equation} 
The above $v$ and $b$ are uniquely determined by the limits
\begin{equation} \label{eq limit defn of v}
    v:= \lim_{t \to +\infty} \dot{x}(t)= \lim_{t \to +\infty} \frac{x(t)}{t}
\end{equation}
and
\begin{equation} \label{eq limit defn of b}
    b:= \begin{cases}
        \lim \limits_{t \to +\infty} \Big[ x(t)- vt- \sum \limits_{j=1}^n c_j(v) t^{1-jp} \Big], & p^{-1} \not\in \N, \\
        \lim \limits_{t \to +\infty} \Big[ x(t)- vt- \sum \limits_{j=1}^n c_j(v) t^{1-jp}- c_{n+1}(v) \log t \Big], & p^{-1} \in \N.
    \end{cases}
\end{equation}
Moreover, the parametrization map $x \mapsto (v,b)$ defined by \eqref{eq limit defn of v} and \eqref{eq limit defn of b} is a bijection from the set of all hyperbolic solutions of \eqref{eq N-body} to $\Y \times \X$.
\end{cor}

\begin{rmk}
In terms of bijectivity, we identify two hyperbolic solutions if there exists $T>0$ such that they agree on $(T,+\infty)$. The coefficients $c_1, \cdots, c_{n+1}$ are identical to those identified in Wu~\cite[Theorem~1]{Wu}.
\end{rmk}

Convergent asymptotic expansions for hyperbolic solutions in the Newtonian case $p=1$ go back to Chazy~\cite{Chazy}. Duignan--Moeckel--Montgomery--Yu~\cite{DMMY} provided a new proof using analytic linearization near infinity. Stable manifold methods have been applied in the parabolic and hyperbolic-parabolic regimes when $p=1$, notably by Saari--Hulkower~\cite{SaariHulkower} and Saari~\cite{Saari}, who studied the manifold structure of initial conditions producing such motions. More recently, Yang-Zhang~\cite{YangZhang} established strong stable-manifold descriptions and uniqueness results for parabolic solutions asymptotic to minimal central configurations when $p=1$.

Our Theorem~\ref{thm rigidity mixed} derives the convergent asymptotic expansions for every profiled expansive solution with an admissible asymptotic pair for $0<p<2$, covering the hyperbolic, parabolic and hyperbolic-parabolic regimes. In the hyperbolic case, Theorem~\ref{thm rigidity hyperbolic} provides a more explicit analytic representation for all $p>0$, extending the Newtonian result in Chazy~\cite{Chazy}; Corollary~\ref{cor hyperbolic} further yields a complete parametrization. Our approach combines stable manifold theory with analytic Poincar\'e--Dulac normal form. We find the argument robust in the sense that it can handle different decay rates, complicated spectral structure and possible resonances arising from the cluster dynamics and the non-Newtonian potential. 

These results are also closely related to Wu~\cite{Wu}, who constructed solutions with prescribed asymptotic data for general homogeneous potentials using a perturbative argument. The two approaches are complementary: Wu~\cite{Wu} constructs solutions with prescribed subleading asymptotic behavior, while we determine the asymptotic structure for arbitrary solutions with admissible profiles. We compare the results in more detail.

In the parabolic case, \cite[Theorem~2]{Wu} identified the upper bound $\frac{p}{p+2}$ for the power modes of the linearized equation and the threshold $\frac{2p-2}{p+2}$ for higher-order corrections. In the hyperbolic-parabolic case, \cite[Theorem~3]{Wu} exhibits two possible correction orders $t^\frac{p}{p+2}$ and $t^{1-p}$ arising from intra-cluster and inter-cluster dynamics, respectively. Our Theorem~\ref{thm rigidity mixed} shows that the corresponding asymptotic structures arise for every solution with an admissible profile and that the full expansion converges. In the hyperbolic case, \cite[Theorem~1]{Wu} established existence and uniqueness for solutions with prescribed asymptotic velocity-offset pair $(v,b)$, which, together with Theorem~\ref{thm rigidity hyperbolic}, yields the bijection assertion in Corollary~\ref{cor hyperbolic}. 

A natural future question is whether the asymptotic rigidity established here, together with the construction in \cite{Wu}, can lead to a complete parametrization of parabolic or hyperbolic-parabolic solutions analogous to the hyperbolic case.

\subsection{Existence and variational properties} \

We next turn to existence and variational properties of profiled expansive solutions. We construct solutions with prescribed initial configurations and asymptotic profiles, and establish their free-time minimizing properties.

\begin{thm}[Minimizing profiled expansive solutions] \label{thm existence mixed} \

Assume $d \in \N_{\ge 2}$ and $p \in (0,2)$. Let $(v,a)$ be an admissible pair and $x_0 \in \X$. Then there exists a solution $x_* \in C^0 \big( [0,+\infty); \X \big) \cap C^2 \big( (0,+\infty); \Y \big)$ of \eqref{eq N-body} with $x_*(0)=x_0$ and
\begin{equation}
    x_*(t)= vt+ at^{\frac{2}{p+2}}+ o(t^{\frac{2}{p+2}}) \quad \text{as} \;\; t \to +\infty.
\end{equation}
Moreover, $x_*$ is a free-time minimizer on $[0, +\infty)$ at energy $|v|_\lambda^2/2$.
\end{thm}

Theorem \ref{thm existence mixed} deals with hyperbolic, parabolic and hyperbolic-parabolic regimes in a unified manner, which requires $p \in (0,2)$. In the hyperbolic case, the conclusion extends to all $p>0$. 

\begin{thm}[Minimizing hyperbolic solutions] \label{thm existence hyperbolic} \

Assume $d \in \N_{\ge 2}$ and $p \ge 2$. Let $v \in \Y$ and $x_0 \in \X$. Then there exists a solution $x_* \in C^0 \big( [0,+\infty); \X \big) \cap C^2 \big( (0,+\infty); \Y \big)$ of \eqref{eq N-body} with $x_*(0)=x_0$ and
\begin{equation}
    x_*(t)= vt+ o(t) \quad \text{as} \;\; t \to +\infty.
\end{equation}
Moreover, $x_*$ is a free-time minimizer on $[0, +\infty)$ at energy $|v|_\lambda^2/2$. 
\end{thm} 

\begin{rmk}
We use the extended-real convention for free-time minimizers, allowing infinite action on compact intervals; see Definition~\ref{def minimizer}.
\end{rmk}

In the Newtonian case $p=1$, the parabolic and hyperbolic existence results corresponding to Theorem~\ref{thm existence mixed} were established by Maderna--Venturelli~\cite{MVparabolic, MVhyperbolic}, respectively. In the hyperbolic-parabolic regime, Burgos \cite{Burgosmixed} constructed solutions with prescribed positive energy and initial configurations, without prescribing the asymptotic profiles. Subsequently, Polimeni--Terracini~\cite{PolimeniTerracini} established the mixed existence results corresponding to Theorem~\ref{thm existence mixed} for $p=1$. For general homogeneous potentials, hyperbolic existence results were obtained by Liu--Yan--Zhou~\cite{LiuYanZhou} and Yu~\cite{Yu}. More recently, Berti--Polimeni--Terracini~\cite{BPT} extended the variational framework in \cite{PolimeniTerracini} to general homogeneous potentials.

Our two theorems above extend and refine the results of \cite[Theorem~1.7, 1.9, 1.10, 1.12]{BPT}. 

In particular, their existence theorem for hyperbolic-parabolic solutions assumes $p \in (\f12,2)$, while Theorem~\ref{thm existence mixed} covers the entire range $p \in (0,2)$. This improvement is based on two ingredients: the construction of a precise enough approximate solution from Wu~\cite[Lemma~4.1]{Wu} (see Lemma~\ref{lem approximate solution}), and a global quadratic lower bound near the reference trajectory in Lemma~\ref{lem lower bound on trajectory}. 

Theorem~\ref{thm existence hyperbolic} addresses the strong-force regime $p \ge 2$ with arbitrary initial configurations, including collisions. This case was not covered by Liu--Yan--Zhou \cite{LiuYanZhou}, which requires $x_0 \in \Y$, or by Yu~\cite{Yu}, which requires $p \in (0,2)$. The key ingredient is a coupled variational formulation that simultaneously handles the singular behavior near the initial collision and the prescribed asymptotics at infinity. 

\begin{rmk}
For $p \ge 2$, any curve with finite Lagrangian action must be collision free, including at the endpoints. Therefore, existing hyperbolic existence results based on finite-action variational results require the initial configuration to be collision-free; see, for example, \cite[Theorem~1.1]{Recent}. We refer to Section~\ref{sec existence hyperbolic} for further discussions.

\end{rmk}

Moreover, the previous rigidity results Theorems~\ref{thm rigidity mixed} and~\ref{thm rigidity hyperbolic} apply to the solutions constructed in Theorem~\ref{thm existence mixed} and~\ref{thm existence hyperbolic}, which cover the remainder estimates and refine the finite-order expansions in \cite{BPT}.
\;\\

The free-time minimizing property is not restricted to the solutions constructed above. In fact, for $p \neq 2$, any profiled expansive solution with an admissible asymptotic pair is eventually a free-time minimizer. 

\begin{thm}[Eventual free-time minimality] \label{thm eventual minimizer} \

Assume $d \in \N_{\ge 1}$ and $p>0$, $p \neq 2$. Let $(v,a)$ be an admissible pair. Assume $x \in C^2 \big( (T,+\infty); \Y \big)$ solves \eqref{eq N-body} on $(T, +\infty)$ for some $T>0$ and
\begin{equation} \label{eq profiled expansive solution}
    x(t) = vt + at^{\frac{2}{p+2}} + o(t^{\frac{2}{p+2}}) \quad \text{as} \;\; t \to +\infty.
\end{equation}
Then there exists $T_0> T$ such that $x$ is a free-time minimizer on $[T_0, +\infty)$ at energy $|v|_\lambda^2/2$.
\end{thm}

Eventual free-time minimality for the Newtonian case $p=1$ was established by Moeckel--Montgomery--S\'anchez Morgado \cite{MMM} for parabolic solutions asymptotic to Lagrange configurations when $N=3$, and by Maderna--Venturelli \cite{MVBusemann} for hyperbolic solutions, both using Busemann functions. Recently, Yang--Zhang~\cite{YangZhang} proved eventual free-time minimality for parabolic solutions asymptotic to minimal central configurations when $p=1$, using stable manifold methods.

Theorem~\ref{thm eventual minimizer} extends these results to all hyperbolic, parabolic and hyperbolic-parabolic regimes, and to general homogeneous potentials. Our argument is variational and relies on the global quadratic lower bound in Lemma~\ref{lem lower bound on trajectory}.
\;\\

The paper is organized as follows. In Section~\ref{sec pre}, we collect several preliminary results, including existence of approximate solutions (Lemma~\ref{lem approximate solution}), and quadratic estimates for the potential (Lemma~\ref{lem Taylor lower bound}--\ref{lem lower bound on trajectory}). Section~\ref{sec asymptotic} is devoted to asymptotic rigidity, where we give the proof of Theorems~\ref{thm rigidity mixed},~\ref{thm rigidity hyperbolic} and Corollary~\ref{cor hyperbolic} following a unified strategy. Finally, in Section~\ref{sec variational}, we prove the existence and eventual minimality results, Theorems~\ref{thm existence mixed},~\ref{thm existence hyperbolic} and~\ref{thm eventual minimizer}. 

\;\\
\textbf{Notations.} In this paper, we let $\N$ denote the set of nonnegative integers. We use the term \textit{polynomial} also for vector-valued polynomials. For two nonnegative functions $f$ and $g$, we write $f \sim g$ if there exists a constant $C>0$ such that $f \le Cg$ and $g \le Cf$. We let $\nabla$ denote the gradient with respect to the mass metric $|\cdot|_\lambda$ and $D$ denote the usual Euclidean differentiation operator. For a matrix $A$, we let $\sigma(A)$ denote the set of its eigenvalues.

\section{Preliminaries} \label{sec pre}

In this section, we collect several preliminary lemmas that will be used in the variational arguments and the rigidity analysis.

We first recall a profiled approximate solution constructed by Wu~\cite{Wu}, which accounts for the slow decay of $U$ at large time when $p$ is small. This is crucial for extending the hyperbolic-parabolic existence result \cite[Theorem~1.12]{BPT} to the full range $p \in (0, 2)$ in Theorem~\ref{thm existence mixed}. 

For $\gamma \in \R$, define the function space 
\begin{equation}
    S_{\gamma}:= \left\{ \sum_{i=1}^m y_i t^{q_i} \log^{k_i} (t) \ \Big| \ m \in \N,\ y_i \in \X,\ q_i \le \gamma,\ k_i \in \N \right\}.
\end{equation}
We note that for $f \in S_\gamma$, we have $f(t)= o(t^{\gamma+})$ as $t \to +\infty$.
The following lemma was proved in \cite[Lemma~4.1]{Wu}. 
\begin{lem} \label{lem approximate solution}
Assume $d \in \N_{\ge 1}$ and $p>0$. Let $(v,a)$ be an admissible pair. Then there exist $f \in S_{1-p}$ such that $r^0(t):= vt+ at^{\frac{2}{p+2}}+ f(t)$ satisfies
\begin{equation} \label{eq approximate solution estimate}
    \ddot{r}^0(t)- \nabla U(r^0(t)) = o(t^{-2}) \quad \text{as} \;\; t \to +\infty.
\end{equation}
\end{lem}

Next, we establish a global lower bound for the potential near any suitable reference trajectory. This helps establish coercivity in the variational arguments. We start with the following pointwise estimate. 

\begin{lem} \label{lem Taylor lower bound}
Assume $d \in \N_{\ge 1}$ and $p>0$. Let $(v,a)$ be an admissible pair and $I \in M_v$. Then for any $\epsilon>0$, there is a neighborhood $V$ of $a_I$ such that for any $b \in V$ and $h \in \X$, we have
\begin{equation} \label{eq Taylor lower bound}
    U_I(b+h)- U_I(b)- (\nabla U_I(b), h)_\lambda \ge -\left( \frac{p}{(p+2)^2}+ \epsilon \right) |h_I|_\lambda^2.
\end{equation}
\end{lem}

\begin{proof}
We first show that for $h \in \X$, we have
\begin{equation} \label{eq Taylor around a_I}
    U_I(a_I+h)- U_I(a_I)- (\nabla U_I(a_I), h)_\lambda \ge -\frac{p}{(p+2)^2} |h_I|_\lambda^2.
\end{equation}
Using Definition~\ref{def admissible pair} and the homogeneity of $U$, we have 
\begin{equation}
    U_I(a_I+h) \ge U_I(a_I) \left( \frac{|a_I|_\lambda}{|a_I+h_I|_\lambda} \right)^p
\end{equation}
and
\begin{equation}
    \nabla U_I(a_I)= -\frac{2p}{(p+2)^2} a_I, \quad U_I(a_I)= \frac{2}{(p+2)^2} |a_I|_\lambda^2.
\end{equation}
By these relations, to prove \eqref{eq Taylor around a_I}, it suffices to prove
\begin{equation}
    \frac{2}{(p+2)^2} \frac{|a_I|_\lambda^{p+2}}{|a_I+h_I|_\lambda^p} - \frac{2}{(p+2)^2} |a_I|_\lambda^2 + \frac{2p}{(p+2)^2} (a_I, h_I)_\lambda \ge -\frac{p}{(p+2)^2} |h_I|_\lambda^2,
\end{equation}
which, by completing the square, is equivalent to
\begin{equation}
    \frac{2}{(p+2)^2} \frac{|a_I|_\lambda^{p+2}}{|a_I+ h_I|_\lambda^p}- \frac{1}{p+2} |a_I|_\lambda^2+ \frac{p}{(p+2)^2} |a_I+ h_I|_\lambda^2 \ge 0.
\end{equation}
Setting $s = |a_I|_\lambda / |a_I+h_I|_\lambda$, then the weighted AM-GM inequality
\begin{equation}
    \frac{2}{p+2} s^p + \frac{p}{p+2} s^{-2} \ge 1
\end{equation}
implies the desired estimate. We thus obtain \eqref{eq Taylor around a_I}. 

Note that by the Taylor formula, \eqref{eq Taylor around a_I} also implies
\begin{equation} \label{eq Hessian lower bound}
    \big( \nabla^2 U_I(a_I)(h), h \big)_\lambda \ge -\frac{2p}{(p+2)^2} |h_I|_\lambda^2,
\end{equation}
where $\nabla^2 U_I$ is with respect to the mass metric $|\cdot|_\lambda$; see \cite[Section~3]{Wu}.

Now go back to the proof of Lemma \ref{lem Taylor lower bound}. We may assume $h=h_I$. Using proof by contradiction, suppose \eqref{eq Taylor lower bound} fails for some $b_n \to a_I$ and $h_n \in \X$ with $h_n= (h_n)_I$. By passing to a subsequence, we may assume $h_n \to h \in \X$ or $|h_n|_\lambda \to +\infty$. Since $U_I>0$ and the right hand side of \eqref{eq Taylor lower bound} is quadratic in $h$, we see $|h_n|_\lambda \to +\infty$ is impossible. If $h_n \to h \neq 0$, then \eqref{eq Taylor around a_I} is violated due to the lower semicontinuity of $U_I$. Therefore, it suffices to prove \eqref{eq Taylor lower bound} when $h$ is close to $0$.

If $b$ is close to $a_I$ and $h$ is close to $0$, then the segment connecting $b$ and $b + h$ is collisionless. By Taylor's formula with the Lagrange remainder, there exists $\theta \in [0,1]$ such that
\begin{equation}
    U_I(b+ h)- U_I(b)- (\nabla U_I(b), h)_\lambda = \frac{1}{2} \big( \nabla^2 U_I(b + \theta h)(h), h \big)_\lambda.
\end{equation}
The continuity of $\nabla^2 U_I$ on the collisionless set and \eqref{eq Hessian lower bound} imply that for $b$ in some neighborhood of $a_I$ and $h$ in some neighborhood of $0$, one has
\begin{equation}
    U_I(b+ h)- U_I(b)- (\nabla U_I(b), h)_\lambda \ge \frac{1}{2} \left( -\frac{2p}{(p+2)^2} - \epsilon \right) |h_I|_\lambda^2.
\end{equation}
This completes the proof of Lemma~\ref{lem Taylor lower bound}.
\end{proof}

Note that during the proof, we obtain a lower bound for the Hessian \eqref{eq Hessian lower bound}. This estimate will be used in the rigidity analysis as well.

The following lemma gives the global quadratic lower bound we need. It applies to arbitrary reference trajectories with admissible asymptotic profiles. 

\begin{lem} \label{lem lower bound on trajectory}
Assume $d \in \N_{\ge 1}$ and $p>0$. Let $(v,a)$ be an admissible pair and 
\begin{equation}
    r \in C^0([0,+\infty); \Y), \quad r(t)= vt+ at^{\frac{2}{p+2}}+ o(t^{\frac{2}{p+2}}) \quad \text{as} \;\; t \to +\infty.
\end{equation}
Then there exists $\kappa \in L^\infty \big( [0,+\infty); \R_{\ge 0} \big)$ such that 
\begin{equation} \label{eq taylor lower bound on trajectory}
    U(r(t)+h)- U(r(t))- (\nabla U(r(t)), h)_\lambda \ge -\kappa(t) |h|_\lambda^2, \;\; \forall t \ge 0, \; h \in \X
\end{equation}
and 
\begin{equation} \label{eq limsup of kappa}
    \limsup_{t \to +\infty} \ t^2 \kappa(t) \le \frac{p}{(p+2)^2}.
\end{equation}
\end{lem}

\begin{proof}
We see that \eqref{eq taylor lower bound on trajectory} holds if we define
\begin{equation}
    \kappa(t)= \sup_{h \in \X, h \neq 0} \frac{\big[ -U(r(t)+h)+ U(r(t))+ (\nabla U(r(t)), h)_\lambda \big]_+}{|h|_\lambda^2}, \;\; \forall t \ge 0.
\end{equation}
We only need to show $\kappa \in L^\infty([0,+\infty); \R_{\ge 0})$ and \eqref{eq limsup of kappa} holds.

Write $r(t)= vt+ at^\frac{2}{p+2}+ f(t)$. Then $f(t)= o(t^{\frac{2}{p+2}})$ as $t \to +\infty$. We decompose
\begin{equation}
    U(x) = \sum_{I \in M_v} U_I(x) + \frac{1}{2} \sum_{I \neq J \in M_v} U_{IJ}(x),
\end{equation}
where for every $I \neq J \in M_v$,
\begin{equation}
    U_{IJ}(x):= \frac{1}{p} \sum_{i \in I, j \in J} \frac{\lambda_i \lambda_j}{|x_i- x_j|^p}.
\end{equation}

For an intra-cluster term $U_I$ with $I \in M_v$, we have
\begin{align}
    U_I(r(t)+h) &= U_I \big( a_I t^{\frac{2}{p+2}}+ f(t)+ h \big) \\
    &= t^{-\frac{2p}{p+2}} U_I \big( a_I+ t^{-\frac{2}{p+2}} f(t)+ t^{-\frac{2}{p+2}} h \big).
\end{align}
Since $f(t)= o(t^{\frac{2}{p+2}})$ as $t \to +\infty$, for any $\epsilon>0$, there exists $T_1>0$ such that for $t \ge T_1$, we can apply Lemma~\ref{lem Taylor lower bound} with $b=a_I+ t^{-\frac{2}{p+2}} f(t)$. We then deduce, for $t \ge T_1$,
\begin{equation} \label{eq Taylor in cluster} \begin{aligned}
    U_I(r(t)+h) &\ge t^{-\frac{2p}{p+2}} \bigg[ U_I \big( a_I+ t^{-\frac{2}{p+2}} f(t) \big) \\
    &\qquad \qquad + \Big( \nabla U_I \big( a_I+ t^{-\frac{2}{p+2}} f(t) \big), t^{-\frac{2}{p+2}} h \Big)_\lambda \\
    &\qquad \qquad -\left( p(p+2)^{-2}+ \epsilon \right) |t^{-\frac{2}{p+2}} h_I|_\lambda^2 \bigg] \\
    &= U_I(r(t))+ \big( \nabla U_I(r(t)), h \big)_\lambda- \Big[ \frac{p}{(p+2)^2}+ \epsilon \Big] t^{-2} |h_I|_\lambda^2.
\end{aligned} \end{equation}

To deal with inter-cluster terms, first note that for $g(y)=|y|^{-p}$, $y \in \R^d \setminus \{0\}$, there exists a constant $C>0$ only depending on $p$ such that 
\begin{equation}
    g(y+u) - g(y) - Dg(y) \cdot u \ge -C |y|^{-p-2} |u|^2, \;\; \forall y \in \R^d \setminus \{0\}, \ u \in \R^d.
\end{equation}
Write $r(t)=(r_1, \cdots r_N)$, $h=(h_1, \cdots h_N)$ and set $y_{ij}= r_i-r_j$, $u_{ij}= h_i-h_j$. In the following, we let $C, c>0$ denote constants independent of $t$ and $h$.  Then for $I \neq J \in M_v$ and $t \ge 0$, 
\begin{equation} \label{eq inter cluster step} \begin{aligned}
    &\quad \ U_{IJ}(r(t)+h)- U_{IJ}(r(t))- (\nabla U_{IJ}(r(t)), h)_\lambda \\
    &= \frac{1}{p} \sum_{i \in I, j \in J} \lambda_i \lambda_j \big( g(y_{ij}+u_{ij})- g(y_{ij})- Dg(y_{ij}) \cdot u_{ij} \big) \\
    &\ge -C \sum_{i \in I, j \in J} |y_{ij}|^{-2-p} |u_{ij}|^2.
\end{aligned} \end{equation}
For $i,j$ in different clusters, there exists $T_2>0$ such that $|y_{ij}| \ge ct$, $\forall t \ge T_2$. Thus, for $t \ge T_2$,
\begin{equation} \label{eq Taylor between clusters}
    U_{IJ}(r(t)+h)- U_{IJ}(r(t))- \big( \nabla U_{IJ}(r(t)), h \big)_\lambda \ge -Ct^{-2-p} \big( |h_I|_\lambda^2+ |h_J|_\lambda^2 \big).
\end{equation}

Combining \eqref{eq Taylor in cluster}, \eqref{eq Taylor between clusters} and the projection identity $\sum_{I \in M_v} |h_I|_\lambda^2 = |h|_\lambda^2$, we deduce that for any $\epsilon>0$, there exists $T_0= \max\{ T_1, T_2\}$ such that
\begin{equation}
    \kappa(t) \le \left( p(p+2)^{-2}+ \epsilon \right)t^{-2} + Ct^{-2-p}, \quad \forall t \ge T_0.
\end{equation}
This then implies \eqref{eq limsup of kappa}. 

Finally, we prove $\kappa \in L^\infty([0,+\infty))$. We only need to consider the compact interval $[0,T_0]$. Similar to \eqref{eq inter cluster step}, for any $t \ge 0$, we have
\begin{align}
    &\quad \ U(r(t)+h)- U(r(t))- (\nabla U(r(t)), h)_\lambda \\
    &= \frac{1}{p} \sum_{i<j} \lambda_i \lambda_j \big( g(y_{ij}+u_{ij})- g(y_{ij})- Dg(y_{ij}) \cdot u_{ij} \big) \\
    &\ge -C \sum_{i<j} |y_{ij}|^{-2-p} |u_{ij}|^2.
\end{align}
Using $r(t) \in \Y$ for $ t \ge 0$, we have $|y_{ij}| \ge c$ for $t \in [0,T_0]$ by compactness. Thus for $t \le T_0$,
\begin{equation}
    U(r(t)+h)- U(r(t))- (\nabla U(r(t)), h)_\lambda \ge -C |h|_\lambda^2, \quad \forall \; h \in \X. 
\end{equation}
This implies $\kappa \in L^\infty([0,+\infty))$ and concludes the proof.
\end{proof}

We conclude this section with some notions and results that are useful for the variational arguments.

We define the \emph{Lagrangian action functional} for a function $x$ on an interval $[t_1,t_2]$ by
\begin{equation} \label{eq lagrangian}
    \L(x; [t_1,t_2])= \int_{t_1}^{t_2} \Big[ \frac{1}{2} |\dot{x}(t)|_\lambda^2+ U \big( x(t) \big) \Big] \dt.
\end{equation}

\begin{defn}[Fixed-time minimizer and free-time minimizer] \label{def minimizer} \

We call $x_*$ a \emph{fixed-time minimizer} on $[t_1, t_2]$ if $x_* \in H^1([t_1,t_2]; \X)$ and $\L(x_*; [t_1, t_2])<+ \infty$, and one has
\begin{equation}
    \L(x; [t_1,t_2]) \ge \L(x_*; [t_1,t_2])
\end{equation}
for every $x \in H^1([t_1, t_2]; \X)$ with $x(t_1)= x_*(t_1)$ and $x(t_2)= x_*(t_2)$.

Let $E \ge 0$. We call $x_*$ a \emph{free-time minimizer} on $[0, +\infty)$ at energy $E$ if 
\begin{equation}
    x_* \in C^0([0,+\infty); \X) \cap H_{\text{loc}}^1((0,+\infty); \X)
\end{equation}
and one has
\begin{equation} \label{eq free time minimizer}
    \L(x; [0,T])+ ET \ge \L(x_*; [t_1,t_2])+ E(t_2-t_1)
\end{equation}
for every $0 \le t_1< t_2< +\infty$, $T>0$ and $x \in H^1([0,T]; \X)$ with $x(0)= x_*(t_1)$ and $x(T)= x_*(t_2)$. 
\end{defn}

Note that we require a fixed-time minimizer to have finite action, but we do not require a free-time minimizer to have finite action, even on compact intervals. In particular, \eqref{eq free time minimizer} should be understood in the extended-real sense.

We recall Marchal's collision-avoidance lemma. It was originally introduced by Marchal \cite{Marchal}, and later refined and generalized by Chenciner \cite{Chenciner} and Ferrario--Terracini \cite{FerrarioTerracini}.

\begin{lem}
Assume $d \in N_{\ge 2}$ and $p>0$. Let $t_1< t_2$ and $x_* \in H^1([t_1, t_2]; \X)$ be a fixed-time minimizer on $[t_1, t_2]$. Then $x_*(t) \in \Y$ for all $t \in (t_1, t_2)$.
\end{lem}

\begin{rmk} \label{rmk Marchal}
By Marchal's lemma and the Euler--Lagrange equation, we deduce that a fixed-time minimizer $x_* \in H^1([t_1,t_2]; \X)$ is a classical solution of \eqref{eq N-body} in $(t_1,t_2)$.
\end{rmk}

\section{Convergent asymptotic expansions} \label{sec asymptotic}

We prove the rigidity results using analytic dynamical system methods. The main strategy is to establish exponential decay in the log-time, reduce the dynamics to a stable manifold, and apply the Poincar\'e--Dulac normal form theory. We first treat the general regime and then the hyperbolic case, where the simpler spectrum admits a more explicit analytic representation.

\subsection{Rigidity of profiled expansive solutions} \label{sec rigidity mixed} \

In this section, we give the proof of Theorem \ref{thm rigidity mixed}. Assume $d \in \N_{\ge 1}$ and $p \in (0,2)$. Let $(v,a)$ be an admissible pair. Assume $x \in C^2((T, +\infty); \Y)$ solves \eqref{eq N-body} and satisfies \eqref{eq expansion profiled expansive}. 

\;\\
\textbf{Step 1: an improved remainder estimate.}

We first show that \eqref{eq expansion profiled expansive} can be improved to assert a genuine lower order remainder.

\begin{lem} \label{lem improved expansion}
There exists $\epsilon>0$ such that
\begin{equation} \label{eq expansion with epsilon}
    x(t)= vt+ at^{\frac{2}{p+2}}+ O(t^{\frac{2}{p+2}- \epsilon}) \quad \text{as} \;\; t \to +\infty.
\end{equation}
\end{lem}

\begin{proof}
Let $h(t)= t^{-\frac{2}{p+2}} (x(t)- vt- at^{\frac{2}{p+2}})$. Then $h(t) \to 0$ as $t \to +\infty$ by \eqref{eq expansion profiled expansive}. We write $U= U_0+ U_1$, where
\begin{equation}
    U_0= \sum_{I \in M_v} U_I, \qquad U_1= \frac{1}{2} \sum_{I \neq J \in M_v} U_{IJ}.
\end{equation}
Since $v_i = v_j$ if $i, j$ are in the same cluster, we have
\begin{align}
    \nabla U(x) &= \nabla U_0 \big( t^{\frac{2}{p+2}}(a+h) \big)+ \nabla U_1 \big( vt+ t^{\frac{2}{p+2}}(a+h) \big) \\
    &= t^{-\frac{2(p+1)}{p+2}} \nabla U_0(a+h)+ t^{-(p+1)} \nabla U_1 \big( v+ t^{-\frac{p}{p+2}} (a+h) \big).
\end{align}
It follows that $h$ satisfies
\begin{equation} \begin{aligned} 
    \ddot{h}= -\frac{4}{p+2} t^{-1} \dot{h} &+ t^{-2} \Big( \nabla U_0(a+h)- \nabla U_0(a)+ \frac{2p}{(p+2)^2} h \Big) \\
    &+ t^{-2-\frac{p(p+1)}{p+2}} \nabla U_1 \big( v+t^{-\frac{p}{p+2}} (a+h) \big).
\end{aligned} \end{equation}
Using the log-time variable $s= \log t$, we have
\begin{equation} \label{eq equation of h in s} \begin{aligned}
    h''= -\frac{2-p}{p+2} h' &+ \nabla U_0(a+h)- \nabla U_0(a)+ \frac{2p}{(p+2)^2} h \\
    &+ e^{-\frac{p(p+1)}{p+2}s} \nabla U_1 \big( v+e^{-\frac{p}{p+2}s} (a+h) \big),
\end{aligned} \end{equation}
where $'$ stands for taking derivatives with respect to $s$. Let \begin{equation}
    W(y)= U_0(y)+ \frac{p}{(p+2)^2} |y|_\lambda^2, \quad g= e^{-\frac{p(p+1)}{p+2}s} \nabla U_1 \big( v+e^{-\frac{p}{p+2}s} (a+h) \big).
\end{equation} 
Note that \eqref{eq Taylor around a_I} and $\nabla U_I(a_I)= -\frac{2p}{(p+2)^2}a_I$ imply
\begin{equation}
    W(a+h) \ge W(a), \quad \nabla W(a)=0.
\end{equation}
We can also simplify \eqref{eq equation of h in s} to
\begin{equation}
    h''= -\frac{2-p}{p+2} h'+ \nabla W(a+h)+ g.
\end{equation}

Let $w= h'$. By integrating $(e^{\frac{2-p}{p+2}s} w)'$, we get
\begin{equation}
    w(s)= e^{-\frac{2-p}{p+2} (s-s_0)} w(s_0)+ \int_{s_0}^s e^{-\frac{2-p}{p+2} (s-\tau)} \big[ \nabla W(a+h(\tau))+ g(\tau) \big] \d\tau.
\end{equation}
Since $h(s) \to 0$ as $s \to +\infty$, we have $\nabla W(a+h(\tau))+ g(\tau) \to 0$ as $\tau \to +\infty$. Then by taking $s_0$ large enough and then $s$ large enough, we see $w(s) \to 0$ as $s \to +\infty$. 

In the following, we let $C, c$ denote positive constants independent of $s$. Define 
\begin{equation}
    \mathcal{E}(s)= \frac{1}{2} |w|_\lambda^2- W(a+h)+ W(a).
\end{equation}
Then $\mathcal{E}(s) \to 0$, $\mathcal{E}(s) \le |w|_\lambda^2/2$ and 
\begin{equation} \label{eq expression of E'}
    \mathcal{E}'(s)= (w, w')_\lambda- (\nabla W(a+h), h')_\lambda= -\frac{2-p}{p+2} |w|_\lambda^2+ (w,g)_\lambda.
\end{equation}
Since $w(s) \to 0$ and $g$ decays exponentially in $s$, we have
\begin{equation}
    \Big| \int_s^\infty (w, g)_\lambda \d\tau \Big|= O(e^{-cs}).
\end{equation}
Note $\mathcal{E}(s) \to 0$. Thus $w \in L^2(s,+\infty)$ for $s$ large enough, and
\begin{equation} \label{eq integral form of E}
    \mathcal{E}(s)= \frac{2-p}{p+2} \int_s^\infty |w|_\lambda^2 \d\tau+ O(e^{-cs}).
\end{equation}
From \eqref{eq expression of E'} and the Cauchy--Schwarz inequality, we also have
\begin{equation}
    \mathcal{E}'(s) \le -c|w|_\lambda^2+ O(e^{-cs}) \le -c\mathcal{E}(s)+ O(e^{-cs}),
\end{equation}
and thus $\mathcal{E}(s) \le Ce^{-cs}$. Then by \eqref{eq integral form of E}, we deduce 
\begin{equation}
    \int_s^\infty |w|_\lambda^2 \d\tau \le Ce^{-cs}.
\end{equation}
Therefore, since $h(s) \to 0$ as $s \to +\infty$, we have
\begin{align}
    |h(s)| &\le \int_s^\infty |w|_\lambda \d\tau= \sum_{k=0}^\infty \int_{s+k}^{s+k+1} |w|_\lambda \d\tau \\
    &\le \sum_{k=0}^{\infty} \Big( \int_{s+k}^{s+k+1} |w|_\lambda^2 \Big)^{\f12} \le \sum_{k=0}^\infty Ce^{-c(s+k)} \le Ce^{-cs}.
\end{align}
This implies \eqref{eq expansion with epsilon} under the original time variable $t=e^s$.
\end{proof}

\noindent
\textbf{Step 2: reduction to the local stable manifold.}

Let $h(t)$ be as in the proof of Lemma \ref{lem improved expansion}. By \eqref{eq equation of h in s}, using the log-time $s= \log t$, we have
\begin{equation} \begin{aligned}
    h'' &= -\frac{2-p}{p+2} h'+ \Big( A+ \frac{2p}{(p+2)^2} \Big) h+ N(h) \\
    &\quad \ + e^{-\frac{p(p+1)}{p+2} s} \nabla U_1 \big( v+e^{-\frac{p}{p+2}s} (a+h) \big),
\end{aligned} \end{equation}
where $A$ is the matrix representation of $\nabla^2 U_0(a)$ with respect to $|\cdot|_\lambda$, and
\begin{equation}
    N(h)= \nabla U_0(a+h)- \nabla U_0(a)- Ah= O(|h|^2).
\end{equation}
Let $w=h'$, $\eta= e^{-\frac{p(p+1)}{p+2} s}$ and $\xi= e^{-\frac{p}{p+2}s}$. Then we obtain the first order system:
\begin{equation} \label{eq first order system} \begin{cases}
    h'= w, \\
    w'= -\frac{2-p}{p+2} w+ (A+ \frac{2p}{(p+2)^2}) h+ N(h)+ \eta \nabla U_1 (v+ \xi(a+h)), \\
    \xi'= -\frac{p}{p+2} \xi, \\
    \eta'= -\frac{p(p+1)}{p+2} \eta.
\end{cases} \end{equation}

Write $u= (h, w, \xi, \eta)^t \in E:= \X \oplus \X \oplus \R \oplus \R$. Then \eqref{eq first order system} has the form 
\begin{equation} \label{eq equation of u}
    u'= Bu+ G(u),
\end{equation}
where the linear part
\begin{equation} 
    B= \begin{pmatrix}
        0 & I & 0 & 0 \\
        A+ \frac{2p}{(p+2)^2} I & -\frac{2-p}{p+2} I & 0 & \nabla U_1(v) \\
        0 & 0 & -\frac{p}{p+2} & 0 \\
        0 & 0 & 0 & -\frac{p(p+1)}{p+2}
    \end{pmatrix}
\end{equation}
and the nonlinear part
\begin{equation}
    G(u)= \big(0, N(h)+ \eta \big[ \nabla U_1 (v+ \xi(a+h))- \nabla U_1(v) \big], 0, 0 \big)^t.
\end{equation}
Since $a_i \neq a_j$ for bodies in the same cluster, and $v_i \neq v_j$ for bodies in different clusters, we deduce that $U_0$, $U_1$ are analytic near $a$, $v$, respectively. Therefore, $G$ is analytic near $u=0$. Moreover, we have $G(0)=0$ and $DG(0)=0$.

By \eqref{eq expansion with epsilon} and the definition of $\eta, \xi$, we see $h, \xi, \eta$ all decay exponentially in $s$. Since $G(u)$ does not depend on $w$ and $G(0)=0$, we have $|G(u)| \le C|(h, \xi, \eta)| \le Ce^{-cs}$ for some $C, c>0$. The equation for $w'$ in \eqref{eq first order system} then gives
\begin{equation}
    w'= -\frac{2-p}{p+2} w+ O(e^{-cs}),
\end{equation}
which yields $w(s)= O(e^{-cs})$ by integrating $(e^{\frac{2-p}{p+2}s} w)'$ and decreasing $c$ if necessary. We deduce
\begin{equation} \label{eq exponential decay of u}
    |u(s)|= O(e^{-cs}) \quad \text{as} \;\; s \to +\infty \;\; \text{for some} \;\; c>0.
\end{equation}

Let $\phi_s$ be the local flow of \eqref{eq equation of u} and $F=\phi_1$. Then $F$ is defined in a neighborhood of $0 \in E$. Moreover, $F$ is analytic with $F(0)=0$ and $DF(0)= e^B$. Note $e^B$ is invertible, so $F$ is a local analytic diffeomorphism. Next, we study the spectrum of $e^B$.

By \cite[Lemma~3.1]{Wu}, $A$ is diagonalizable. Thus there are eigenvectors $b_1, \cdots, b_{dN}$ of $A$ forming a basis of $\X$, with eigenvalues $c_1, \cdots, c_{dN}$, respectively. Recalling \eqref{eq Hessian lower bound} or using \cite[Lemma~3.1]{Wu}, we have 
\begin{equation}
    c_i \ge -\frac{2p}{(p+2)^2}, \quad \forall 1 \le i \le dN.
\end{equation}
In particular, $c_i> -\f14$, so the quadratic equation $\lambda^2- \lambda= c_i$ has two distinct real roots. We denote them by $p_i, q_i$ with $p_i< q_i$. Using $c_i \ge -\frac{2p}{(p+2)^2}$, we have 
\begin{equation}
    p_i \le \frac{p}{p+2}, \quad q_i \ge \frac{2}{p+2}.
\end{equation} Let $\mu_i^-= p_i- \frac{2}{p+2}$, $\mu_i^+= q_i- \frac{2}{p+2}$ and
\begin{equation} \label{eq defn of B0}
    B_0= \begin{pmatrix}
        0 & I  \\
        A+ \frac{2p}{(p+2)^2} I & -\frac{2-p}{p+2} I
    \end{pmatrix}.
\end{equation}
Then $\mu_i^- \le -\frac{2-p}{p+2} <0$, $\mu_i^+ \ge 0$ and
\begin{equation}
    B_0 \begin{pmatrix}
        b_i \\
        \mu_i^\pm b_i
    \end{pmatrix} = \mu_i^\pm\begin{pmatrix}
        b_i \\
        \mu_i^\pm b_i
    \end{pmatrix}.
\end{equation}
Thus $B_0$ is diagonalizable and its eigenvalues are $\{\mu_i^\pm\}_{i=1}^{dN}$. Since $B$ is block upper-triangular, 
\begin{equation}
    \sigma(B)= \bigcup_{i=1}^{dN}  \{\mu_i^-, \mu_i^+\} \cup \Big\{ -\frac{p}{p+2}, -\frac{p(p+1)}{p+2}  \Big\},
\end{equation}
so $B$ has $dN+2$ negative eigenvalues and $dN$ nonnegative eigenvalues, counted with multiplicity. Define
\begin{equation}
    E^-= \bigoplus_{\lambda \in \sigma(B), \lambda<0} \ker \big( (B- \lambda)^{2dN+2} \big), \;\; E^+= \bigoplus_{\lambda \in \sigma(B), \lambda \ge 0} \ker \big( (B- \lambda)^{2dN+2} \big).
\end{equation}
Then $E^-, E^+$ are $B$-invariant, $\dim E^-= dN+2$, $\dim E^+= dN$ and $E= E^- \oplus E^+$. We point out that $B|_{E^-}$ may not be diagonalizable because of $\nabla U_1(v)$, but $B|_{E^+}$ is diagonalizable.

We have $\sigma(e^B | _{E^-}) \subset \{ z \ | \ |z|<1 \}$ and $\sigma(e^B | _{E^+}) \subset \{ z \ | \ |z| \ge 1 \}$. Moreover, since 
\begin{equation}
    \sup \{ |\lambda| \ | \ \lambda \in \sigma(e^B |_{E^-}) \} \le \max \big\{ e^{-\frac{2-p}{p+2}}, e^{-\frac{p}{p+2}}, e^{-\frac{p(p+1)}{p+2}} \big\} <1,
\end{equation}
there exists a positive integer $L$ such that
\begin{equation}
    \big( \sigma(e^B | _{E^-}) \big)^{L+1} \sigma (e^{-B}) \subset \{ z \ | \ |z|<1 \}.
\end{equation}
For any $2 \le i \le L$, we have
\begin{equation}
    \big( \sigma(e^B | _{E^-}) \big)^i \cap \sigma (e^B |_{E^+}) = \varnothing.
\end{equation}
Since $F$ is a local diffeomorphism, all the assumptions of \cite[Theorem~1.1]{LocalInvManifold} are satisfied. We conclude that there exists a local $F$-invariant analytic manifold $W^- \subset E$ tangent to $E^-$ at $0$. By \cite[Remark~6]{LocalInvManifold}, $W^-$ is also locally invariant under the flow of \eqref{eq equation of u}. 

Since $W^-$ is tangent to $E^-$ at $0$, there exist a neighborhood $V$ of $0$ in $E^-$, a neighborhood $\tilde{V}$ of $0$ in $W^-$ and an analytic function $\psi: V \to E^+$ such that $\psi(0)=0$, $D\psi(0)=0$ and 
\begin{equation}
    \tilde{V}= \{ y+\psi(y) \ | \ y \in V \}.
\end{equation}
Let $\pi_-, \pi_+$ be the projections onto $E^-, E^+$, respectively. Take $y \in V$ and let $\tilde{u}$ be the solution to \eqref{eq equation of u} with $\tilde{u}(0)= y+ \psi(y)$. By the invariance of $W^-$ under the flow of \eqref{eq equation of u}, we have $\tilde{u} \in W^-$ locally. Thus $\pi_+ \tilde{u}= \psi(\pi_- \tilde{u})$. Differentiating both sides in $s$ and using \eqref{eq equation of u}, we obtain
\begin{equation} \label{eq local invariance of W^-} \begin{aligned}
    &\quad \ B|_{E^+} \psi(y)+ \pi_+ G(y+\psi(y)) \\
    &= D\psi(y) \big( (B|_{E^-}) y+ \pi_- G(y+\psi(y)) \big), \quad \forall y \in V.
\end{aligned} \end{equation}

\begin{lem} \label{lem eventually stable}
For $s$ large enough, we have $u(s) \in W^-$.
\end{lem}
\begin{proof}
Let $y(s)= \pi_-(u(s))$ and $\delta(s)= \pi_+(u(s))- \psi(y(s))$. Then 
\begin{equation}
    u(s) \in W^- \iff \delta(s)=0.
\end{equation} 
By \eqref{eq equation of u} and invariance of $E^\pm$ under $B$, we have
\begin{equation}
    y'= (B|_{E^-}) y+ \pi_-G(u), \quad
    \pi_+(u(s))'= (B|_{E^+}) \pi_+(u(s)) + \pi_+G(u).
\end{equation}
Then we compute
\begin{equation}
    \delta'= (B|_{E^+}) \pi_+(u(s))+ \pi_+G(u)- D\psi(y) \big( (B|_{E^-}) y+ \pi_-G(u) \big).
\end{equation}
By \eqref{eq local invariance of W^-}, for $s$ large enough, we further have
\begin{equation} \label{eq tem eq of delta}
    \delta'= (B|_{E^+})\delta+ (\pi_+- D\psi(y) \pi_-) \big( G(u)- G(u-\delta) \big).
\end{equation}
From $u(s) \to 0$, we know $\delta(s) \to 0$. Using \eqref{eq exponential decay of u} and $DG(0)=0$, for large $s$, we have 
\begin{equation}
    \Big| (\pi_+- D\psi(y) \pi_-) \big( G(u)- G(u-\delta) \big) \Big| \le C|u||\delta| \le Ce^{-cs} |\delta|.
\end{equation}
Note that $B|_{E^+}$ is diagonalizable and its eigenvalues are nonnegative. There exists an inner product $\langle \cdot, \cdot \rangle$ on $E_+$ such that $\langle B|_{E^+} z, z \rangle \ge 0$ for all $z \in E^+$. By \eqref{eq tem eq of delta}, for large $s$, we have
\begin{equation}
    \frac{\d}{\d s} \langle \delta, \delta \rangle= 2\langle \delta', \delta \rangle \ge -Ce^{-cs} \langle \delta, \delta \rangle.
\end{equation}
Multiplying by an integration factor, for $S_0<S_1$ large enough, we get  
\begin{equation}
    \langle \delta(S_1), \delta(S_1) \rangle \exp \big( \int_{S_0}^{S_1} Ce^{-c\tau} \d\tau \big) \ge \langle \delta(S_0), \delta(S_0) \rangle.
\end{equation}  
Fix $S_0$ and take $S_1 \to +\infty$. Since $e^{-c\tau}$ is integrable at infinity and $\delta(s) \to 0$ as $s \to +\infty$, we deduce $\delta(S_0)=0$. Thus $u(S_0) \in W^-$. By invariance of $W^-$ under the flow of \eqref{eq equation of u}, we get the result of Lemma~\ref{lem eventually stable}.
\end{proof}

\noindent
\textbf{Step 3: normal form and convergent expansion.}

By Lemma~\ref{lem eventually stable}, we have $u(s)= y(s)+ \psi(y(s))$ and it suffices to study the differential equation
\begin{equation} \label{eq equation of y}
    y'= (B|_{E^-})y+ \pi_-G(y+ \psi(y)), \quad y(s) \in V \subset E^-, \quad s \gg 1.
\end{equation}
Let $k:=\dim E^-= dN+2$ and let $\nu_1, \cdots, \nu_k$ be the eigenvalues of $B|_{E^-}$, counted with multiplicity. By the definition of $E^-$, all these eigenvalues are negative. Thus $\nu=(\nu_1, \cdots, \nu_k)$ lies in the Poincar\'e domain. By \cite[Chapter~I, Theorem~5.5]{PoincareDulac}, there is an analytic local diffeomorphism $y= \phi(z)$ such that \eqref{eq equation of y} is transformed to the Poincar\'e-Dulac normal form
\begin{equation} \label{eq equation of z}
    z'= Jz+ P(z),
\end{equation}
where $J$ is the upper-triangular Jordan normal form of $B|_{E^-}$ and $P$ is a polynomial consisting of resonant monomials of degree at least two. More precisely, if ordering the diagonal entries of $J$ exactly as $\nu_1, \cdots, \nu_k$ and writing $P(z)= (P_1(z), \cdots, P_k(z))$, then for each $j$, 
\begin{equation}
    P_j(z)= \sum_{m \in \N^k, |m| \ge 2, m \cdot \nu = \nu_j} a_{j,m} z_1^{m_1} z_2^{m_2} \cdots z_k^{m_k}.
\end{equation}
Here, resonant corresponds to the condition $m \cdot \nu= \nu_j$. 

\begin{lem} \label{lem normal form solution}
If $z(s)$ solves \eqref{eq equation of z}, then there are polynomials $R_1, \cdots, R_k$ such that 
\begin{equation} \label{eq polynomial expression of z}
    z_j(s)= e^{\nu_j s} R_j(s), \quad \forall 1 \le j \le k.
\end{equation}
\end{lem}

\begin{proof}
We will prove the result by induction in $k$. The case when $k=1$ is clear because the equation is linear. Now, assume $k>1$ and the result holds for any dimension less than $k$.

If all $\nu_j$ are equal, then \eqref{eq equation of z} is linear because $m \cdot \nu= \nu_j$ does not admit any $m$ with $|m| \ge 2$. The conclusion is thus clear. Otherwise, assume $\nu_1 \ge \cdots \ge \nu_k$ and take $k'$ to be the largest index so that $\nu_{k'}> \nu_k$. We have $1 \le k' < k$. Write $J$ in the block form
\begin{equation}
    J= \begin{pmatrix}
        J_0 & 0 \\
        0 & \tilde{J}
    \end{pmatrix},
\end{equation}
where $J_0$ is a $k'$ by $k'$ matrix and $\tilde{J}$ is a $(k-k')$ by $(k-k')$ matrix. Then the diagonal entries of $J_0$ are all greater than $\nu_k$ and the diagonal entries of $\tilde{J}$ are all equal to $\nu_k$. Consider the resonant condition $m \cdot \nu= \nu_j$. Since all $\nu_1, \cdots, \nu_k$ are negative, if $m_i>0$ for some $i$, then we must have $\nu_i>\nu_j$. Therefore, for any $1 \le j \le k$, $P_j(z)$ only depends on $z_1, \cdots, z_{k'}$. In particular, the equation for $z_1, \cdots, z_{k'}$ is a sub-problem. Since $k'<k$, by the induction hypothesis, there are polynomials $R_1, \cdots, R_{k'}$ such that
\begin{equation} \label{eq formula of z by induction}
    z_j(s)= e^{\nu_j s} R_j(s), \quad \forall 1 \le j \le k'.    
\end{equation}
For $k'+1 \le j \le k$, since $P_j(z)$ also only depends on $z_1, \cdots, z_{k'}$, by \eqref{eq formula of z by induction} and the resonant condition that $z^m$ appears in $P_j(z)$ only when $m \cdot \nu= \nu_j= \nu_k$, we have
\begin{equation} 
    P_j(z)= e^{\nu_k s} \tilde{P}_j(s), \quad \forall k'+1 \le j \le k,    
\end{equation}
where $\tilde{P}_j$ is a polynomial. Then $\tilde{z}= (z_{k'+1}, \cdots, z_k)$ satisfies the differential equation
\begin{equation}
    \tilde{z}'= \tilde{J} \tilde{z}+ e^{\nu_k s} \tilde{P}(s).
\end{equation}
This is an inhomogeneous linear equation, and its solution is given by
\begin{equation}
    \tilde{z}(s)= e^{(s-s_0)\tilde{J}} \tilde{z}(s_0)+ \int_{s_0}^s e^{(s-\tau) \tilde{J}} e^{\nu_k \tau} \tilde{P}(\tau) \d\tau.
\end{equation}
Recall $\tilde{J}$ consists of the Jordan blocks associated to $\nu_k$, so $\tilde{J}- \nu_k$ is nilpotent, and thus $e^{s (\tilde{J}- \nu_k)}$ is a polynomial in $s$. We write
\begin{equation}
    e^{-\nu_k s} \tilde{z}(s)= e^{s (\tilde{J}- \nu_k)} e^{-s_0 \tilde{J}} \tilde{z}(s_0)+ e^{s (\tilde{J}- \nu_k)} \int_{s_0}^s  e^{-\tau (\tilde{J}- \nu_k)} \tilde{P}(\tau) \d\tau.
\end{equation}
Then the right hand side is a polynomial in $s$. This completes the proof.
\end{proof}

Now we have $u(s)= \phi(z(s))+ \psi(\phi(z(s)))$ and $z$ satisfies \eqref{eq polynomial expression of z}. Let $\Theta(z)= \phi(z)+ \psi(\phi(z))$. By the analyticity of $\phi, \psi$ and $\phi(0)= \psi(0)= 0$, we see that $\Theta$ is analytic near $0$ and $\Theta(0)=0$. Since $z(s) \to 0$ as $s \to +\infty$, for $s$ large enough, the Taylor series gives
\begin{equation}
    u(s)= \sum_{m \in \N^k \setminus \{0\}} a_m z(s)^m \quad \text{for some} \;\; a_m \in E.
\end{equation}
By \eqref{eq polynomial expression of z}, for large $s$, we further have
\begin{equation} 
    u(s)= \sum_{m \in \N^k \setminus \{0\}} e^{(m \cdot \nu) s} a_m R_1(s)^{m_1} \cdots R_k(s)^{m_k}.
\end{equation}
Let $Y=(h,w)^t$ and restrict the expansion to $Y$. By collecting terms with the same exponent, there are real numbers $\beta_1 > \beta_2> \cdots$ and nonzero polynomials $Q_1, Q_2, \cdots$ such that
\begin{equation} \label{eq expansion of Y}
    Y(s)= \sum_i e^{\beta_i s} Q_i(s), \quad \forall s \gg 1,
\end{equation}
where the sum can be empty, finite, or countably infinite. By analyticity of $\phi$ and $\psi$, the series converges absolutely and uniformly for $s$ large enough.

\;\\
\textbf{Step 4: properties of the exponents.} 

Assume the sum in \eqref{eq expansion of Y} is nonempty. Since all $\nu_j<0$, we have $\beta_i<0$. We also have the following more precise conclusions.

\begin{lem} \label{lem tail bound}
For each $i$, we have

(1)
\begin{equation}
    \beta_i \le \max\{ -\frac{2-p}{p+2}, -\frac{p(p+1)}{p+2} \}. 
\end{equation}

(2)
\begin{equation}
    \begin{cases}
        \beta_i \ge -\frac{p(p+1)}{p+2} \\
        \beta_i> -\frac{2(2-p)}{p+2}
    \end{cases} \implies Q_i \text{ is a constant}.
\end{equation}
\end{lem} 

\begin{proof}
Consider the equations for $h$ and $w$ in \eqref{eq first order system}
\begin{equation}
    \begin{cases}
        h'= w, \\
        w'= -\frac{2-p}{p+2} w+ (A+ \frac{2p}{(p+2)^2}) h+ N(h)+ \eta \nabla U_1 (v+ \xi(a+h)).
    \end{cases}
\end{equation}

(1) Assume $\beta_1 > -\frac{p(p+1)}{p+2}$. Then $e^{\beta_1 s}$ is the leading exponential term, and both $N(h)$ and $\eta \nabla U_1 (v+ \xi(a+h))$ are strictly of lower orders. By comparing the leading coefficients, we get 
\begin{equation}
    Q_1'= (B_0- \beta_1) Q_1,
\end{equation}
with $B_0$ defined in \eqref{eq defn of B0}. Thus 
\begin{equation}
    Q_1(s)= e^{(B_0- \beta_1)(s-s_0)} Q_1(s_0) \quad \text{for} \;\; s> s_0 \gg 1.
\end{equation}
Since $Q_1$ is a polynomial, there is a positive integer $D$ such that $(B_0- \beta_1)^D Q_1(s_0)= 0$. Then $B_0$ being diagonalizable implies 
\begin{equation}
    (B_0- \beta_1) Q_1(s_0)= 0
\end{equation} 
and thus $Q_1$ is a constant. Note also $Q_1 \neq 0$, so $\beta_1 \in \sigma(B_0)$. We deduce $\beta_1 \le -\frac{2-p}{p+2}$ because 
\begin{equation}
    \sigma(B_0) \cap (-\infty, 0) \subset (-\infty, -\frac{2-p}{p+2}]. 
\end{equation}
This proves (1) since $\beta_i$ is decreasing.

(2) Assume $\beta_i \ge -\frac{p(p+1)}{p+2}$ and $\beta_i> -\frac{2(2-p)}{p+2}$. Then $\beta_i> 2\beta_1$ by (1). Since $N(h)= O(|h|^2)$, we see $N(h)$ is strictly of lower order than $e^{\beta_i s} Q_i(s)$. 

If $\beta_i > -\frac{p(p+1)}{p+2}$, then $\eta \nabla U_1(v+ \xi(a+h))$ is also of lower order than $e^{\beta_i s} Q_i(s)$. By comparing the coefficients of $e^{\beta_i s}$, we get $Q_i'= (B_0- \beta_i) Q_i$. We can then argue as in (1) to deduce that $Q_i$ is a constant. 

If $\beta_i= -\frac{p(p+1)}{p+2}$, then comparing the coefficients of $e^{\beta_i s}$ yields
\begin{equation} \label{eq coefficient equation with resonance}
    Q_i'= (B_0- \beta_i) Q_i+ (0, \nabla U_1(v))^t.
\end{equation}
Note that $\nabla U_1(v)$ is constant on each cluster, that is, $(\nabla U_1(v))_i= (\nabla U_1(v))_j$ for $i \sim_v j$. By the translation invariance of $\nabla U_0$, we have $A\nabla U_1(v)=0$. Besides, $\beta_i= -\frac{p(p+1)}{p+2}$ and $\beta_i> -\frac{2(2-p)}{p+2}$ implies $p<1$. In particular, $p \neq 1$. By direct calculation,
\begin{equation}
    Q_*= \frac{1}{p(p-1)} (\nabla U_1(v), \beta_i \nabla U_1(v))^t
\end{equation}
is a concrete solution of \eqref{eq coefficient equation with resonance}. Thus
\begin{equation}
    (Q_i- Q_*)'= (B_0- \beta_i)(Q_i- Q_*).
\end{equation}
We deduce $Q_i- Q_*$ is a constant by the argument in (1), so $Q_i$ is a constant as well.
\end{proof}

Now we can derive Theorem~\ref{thm rigidity mixed}.

\begin{proof}[Proof of Theorem~\ref{thm rigidity mixed}]

Under the original time $t= e^s$, \eqref{eq expansion of Y} implies
\begin{equation}
    Y(t)= \sum_i t^{\beta_i} Q_i(\log t), \quad \forall t \gg 1.
\end{equation}
Write $Q_i= (H_i, W_i)^t$, with both $H_i$ and $W_i$ taking values in $\X$. Then
\begin{equation}
    h(t)= \sum_i t^{\beta_i} H_i(\log t), \quad w(t)= \sum_i t^{\beta_i} W_i(\log t), \quad \forall t \gg 1.
\end{equation}
By \eqref{eq first order system}, we have $H_i'+ \beta_i H_i= W_i$. Thus $Q_i \neq 0$ implies $H_i \neq 0$. Using 
\begin{equation}
    x(t)= vt+ at^{\frac{2}{p+2}}+ t^{\frac{2}{p+2}}h(t),
\end{equation} 
we get 
\begin{equation} 
    x(t) = vt + at^{\frac{2}{p+2}}+ \sum_i t^{\alpha_i} H_i(\log t), \quad \forall t \gg 1,
\end{equation}
where $\alpha_i= \beta_i+ \frac{2}{p+2}$ and $H_i$ is a nonzero polynomial. Since the expansion of $Y$ comes from the Taylor series of an analytic function, the series converges absolutely and uniformly for $t$ large enough. Note that all $\beta_i$ belong to the set
\begin{equation}
    \big\{ m \cdot \nu \ | \ m \in \N^k \setminus \{0\} \big\},
\end{equation}
and $\nu$ consists of the negative eigenvalues of $B$. Since $\sigma(B)$ is determined by $a$ and the cluster partition $M_v$, each $\alpha_i$ belongs to a discrete subset of $\R$ that depends only on $a$ and $M_v$. By Lemma~\ref{lem tail bound}, we also obtain \eqref{eq mixed exponent bound} and \eqref{eq mixed no resonant}.

Now, to prove Theorem~\ref{thm rigidity mixed}, it remains to prove the strengthened estimates in the parabolic case. This is a natural improvement since the inter-cluster term vanishes. 

We have $U_1=0$ because $|M_v|=1$. Therefore, the comparison between $\beta_i$ and $-\frac{p(p+1)}{p+2}$ in Lemma~\ref{lem tail bound} is not needed. We then deduce 
\begin{equation}
    \beta_i \le -\frac{2-p}{p+2}
\end{equation} 
and
\begin{equation}
    \beta_i> -\frac{2(2-p)}{p+2} \implies Q_i \; \text{is a constant}.
\end{equation} 
Using $\alpha_i= \beta_i+ \frac{2}{p+2}$, we obtain the desired estimates in the parabolic case.
\end{proof}

\subsection{Rigidity and parametrization of hyperbolic solutions} \

We now consider the hyperbolic case. The strategy is similar to that of Theorem~\ref{thm rigidity mixed}. The condition $v \in \Y$ simplifies the calculation and allows for more explicit expressions. 

Let $v \in \Y$. Assume $x \in C^2 \big( (T,+\infty); \Y \big)$ solves \eqref{eq N-body} and satisfies \eqref{eq expansion hyperbolic}. We calculate the equation for $h(t)= t^{-1} (x(t)-vt)$. Note $h(t)= o(1)$ as $t \to +\infty$ since \eqref{eq expansion hyperbolic}. By \eqref{eq N-body}, we have
\begin{equation}
    \ddot{h}= -2t^{-1} \dot{h}+ t^{-2-p} \nabla U(v+h).
\end{equation}
Using the log-time variable $s=\log t$, we get
\begin{equation}
    h''= -h'+ e^{-ps} \nabla U(v+h),
\end{equation}
where $'$ stands for taking derivatives with respect to $s$.

Let $w=h'$ and $\eta= e^{-ps}$. Then we obtain the first order system:
\begin{equation} \label{eq first order system hyperbolic}
    \begin{cases}
        h'= w \\
        w'= -w+ \eta \nabla U(v+h), \\
        \eta' = -p\eta.
    \end{cases}
\end{equation}
Write $u=(h,w,\eta)^t \in E:= \X \oplus \X \oplus \R$. Then
\begin{equation} \label{eq equation of u hyperbolic}
    u'= Bu+ G(u),
\end{equation}
where the linear part
\begin{equation}
    B= \begin{pmatrix}
        0 & I & 0 \\
        0 & -I & \nabla U(v) \\
        0 & 0 & -p
    \end{pmatrix}
\end{equation}
and the nonlinear part
\begin{equation}
    G(u)= \big( 0, \eta \big[ \nabla U(v+h)- \nabla U(v) \big], 0 \big)^t.
\end{equation}
The function $G$ is analytic near $0$ with $G(0)= 0$ and $DG(0)=0$. 

Next, we show that $u(s)$ decays exponentially in $s$. Using \eqref{eq first order system hyperbolic}, we have 
\begin{equation}
    w(s)= e^{-(s-s_0)} w(s_0)+ e^{-s} \int_{s_0}^s e^{(1-p)\tau} \nabla U(v+ h(\tau)) \d\tau.
\end{equation}
Since $h(s) \to 0$ as $s \to +\infty$, we have $|\nabla U(v+ h(\tau))| \le C$ for $\tau$ large enough. Thus we obtain
\begin{equation}
    |w(s)| \le Cse^{-\min\{1,p\}s} \quad \text{for} \; s \gg 1.
\end{equation}
Using this and $h'=w$, we see $h$ also decays exponentially. We conclude
\begin{equation} \label{eq exponential decay of u hyperbolic}
    |u(s)|= O(e^{-cs}) \quad \text{as} \;\; s \to +\infty \;\; \text{for some} \;\; c>0.
\end{equation}

Then we construct the local stable manifold $W^-$. We have $\sigma(B)= \{0, -1, -p\}$. Let 
\begin{equation}
    E^-= \bigoplus_{\lambda \in \sigma(B), \lambda<0} \ker \big( (B- \lambda)^{2dN+2} \big), \;\; E^+= \bigoplus_{\lambda \in \sigma(B), \lambda \ge 0} \ker \big( (B- \lambda)^{2dN+2} \big).
\end{equation}
By the same argument as in Section~\ref{sec rigidity mixed}, there exists a local analytic manifold $W^-$ tangent to $E^-$ at $0$. Moreover, $W^-$ is locally invariant under the flow of \eqref{eq equation of u hyperbolic}. Using \eqref{eq exponential decay of u hyperbolic}, the same argument in Lemma~\ref{lem eventually stable} implies $u(s) \in W^-$ for $s$ large enough. 

Therefore, there exists an analytic function $\psi$ such that $u(s)= y(s)+ \psi(y(s))$ and \eqref{eq equation of u hyperbolic} is equivalent to 
\begin{equation}
    y'= (B|_{E^-}) y+ \pi_- G(y+ \psi(y)), \quad \forall s \gg 1.
\end{equation}
Here $k:= \dim E^-= dN+1$ and 
\begin{equation}
    \nu_1= \cdots =\nu_{k-1}= -1, \quad \nu_k=-p
\end{equation} 
are all eigenvalues of $B|_{E^-}$ with multiplicity. We can further find an analytic local diffeomorphism $y= \phi(z)$ reducing the equation to the Poincar\'e-Dulac normal form
\begin{equation}
    z'= Jz+P(z),
\end{equation}
where $J$ is the Jordan normal form of $B|_{E^-}$ and $P$ is a polynomial consisting of resonant monomials of degree at least two.

Using the expression of $B$, we have more precise results than Lemma~\ref{lem normal form solution}. Recall $\nu_i$ is either $-1$ or $-p$.

\begin{lem} \label{lem generators}
If $p^{-1} \not\in \N$, then each $z_i$ is a linear combination of $e^{-s}$ and $e^{-ps}$. If $p ^{-1} \in \N$, then each $z_i$ is a linear combination of $e^{-s}$, $e^{-ps}$ and $e^{-s}s$.
\end{lem}

\begin{proof}
First, note $\operatorname{rank} (B+I)|_{E^-} \le 1$, so $J$ is diagonal if $p \neq 1$, and $(J+I)^2=0$ if $p=1$.

If $p=1$, then all eigenvalues are equal, and thus there is no resonant monomial of degree at least two. This means $P(z)=0$, so the equation is simply $z'= Jz$. The solution is
\begin{equation}
    z(s)= e^{-s} e^{s(J+1)} e^{-s_0 J} z(s_0).
\end{equation}
Since $(J+I)^2=0$, the function $e^{s(J+I)}$ is a polynomial of degree at most one. We obtain
\begin{equation}
    z_i(s)= e^{-s}(A_i+ B_i s), \quad \forall 1 \le i \le k,
\end{equation}
where $A_i, B_i$ are constants.

Next assume $p \neq 1$, which implies $J= \text{diag} (\nu_1, \cdots, \nu_k)$. If $p \not \in \N$ and $p^{-1} \not\in \N$, then there is still no resonant monomial of degree at least two. Thus $P(z)=0$ and we directly obtain
\begin{equation}
    z_i(s)= A_i e^{\nu_i s}, \quad \forall 1 \le i \le k.
\end{equation}

If $p^{-1} \in \N_{>1}$, then the resonant condition $m \cdot \nu= \nu_i$ forces $1 \le i \le k-1$, $m_1= \cdots m_{k-1}=0$ and $m_k= p^{-1}$. Therefore, $P_k= 0$ and $P_i(z)= C_i z_k^{1/p}$ for $1 \le i \le k-1$. The equation becomes
\begin{equation}
    \begin{cases} 
        z_k'= -pz_k, \\
        z_i'= -z_i+ C_iz_k^{1/p}, \quad 1 \le i \le k-1.
    \end{cases}
\end{equation}
We can solve it to get 
\begin{equation}
    \begin{cases}
        z_k(s)= A_k e^{-ps}, \\
        z_i(s)= e^{-s} (A_i+ B_is), \quad 1 \le i \le k-1.
    \end{cases}
\end{equation}

If $p \in \N_{>1}$, then the resonant condition $m \cdot \nu= \nu_i$ forces $i=k$, $m_k=0$ and $\sum_{i=1}^{k-1} m_i= p$. Therefore, $P_i(z)=0$ for $1 \le i \le k-1$ and $P_k$ only depends on $z_1, \cdots, z_{k-1}$ and is homogeneous of degree $p$. We further claim that $P_k(z)=0$. Assuming this claim, the equation is linear and we obtain
\begin{equation}
    z_i(s)= A_i e^{\nu_i s}, \quad \forall 1 \le i \le k.
\end{equation}

Combining the above cases, we get the conclusion.

Finally, we supply the proof of $P_k=0$ when $p \in \N_{>1}$. By $u= y+ \psi(y)$ and $y= \phi(z)$, we have $\eta= \theta(z)$, where $\theta$ is analytic near $0$. Then $\eta'= -p\eta$ becomes 
\begin{equation} \label{eq tem identity of z}
    D\theta(z) Jz+ D\theta(z) P(z)= -p\theta(z).
\end{equation} 
We expand $\theta(z)$ into power series and compare the terms in \eqref{eq tem identity of z} depending only on $z_1, \cdots, z_{k-1}$ and homogeneous of degree $p$. By $\nu_1, \cdots, \nu_{k-1}=-1$ and homogeneity, these terms vanish in $D\theta(z) Jz+ p\theta(z)$; moreover, such terms in $D\theta(z) P(z)$ equal $\partial_k \theta(0) P_k$. Thus $\partial_k \theta(0) P_k=0$. Since both $y \mapsto y+\psi(y)$ and $\phi$ are local diffeomorphisms, and $\eta$ has nonzero projection on $E^-$, we have $D\theta(0) \neq 0$. For $1 \le i \le k-1$, by comparing the $z_i$ term in \eqref{eq tem identity of z} and using $p \neq 1$, we get $\partial_i \theta(0)=0$. Thus $\partial_k \theta(0) \neq 0$ and we obtain $P_k=0$.
\end{proof}

Note that if $p^{-1} \in \N$, then $e^{-s}$ is an integer power of $e^{-ps}$. We deduce that $z$ is a polynomial function of $(e^{-ps}, e^{-s})$ if $p^{-1} \not \in \N$, and a polynomial function of $(e^{-ps}, e^{-s}s)$ if $p^{-1} \in \N$.

Theorem~\ref{thm rigidity hyperbolic} now follows immediately.

\begin{proof}[Proof of Theorem~\ref{thm rigidity hyperbolic}]

Restricting $u= \phi(z)+ \psi(\phi(z))$ to $h$, we deduce that there exists a function $\Phi(\cdot, \cdot)$ that is analytic near $(0,0)$ such that
\begin{equation} \label{eq analytic expansion of h}
    h(s)= \begin{cases}
        \Phi(e^{-ps}, e^{-s}), & p^{-1} \not\in \N, \\
        \Phi(e^{-ps}, e^{-s}s), & p^{-1} \in \N
    \end{cases} \quad \text{for} \;\; s \gg 1.
\end{equation}
By $h(s) \to 0$ as $s \to +\infty$, we also deduce $\Phi(0,0)= 0$. Then, by reversing to the original time variable $t= e^s$ and using $x(t)= vt+ th(t)$, we obtain
\begin{equation} 
    x(t)= \begin{cases}
        vt+ t\Phi(t^{-p}, t^{-1}), & p^{-1} \not\in \N, \\
        vt+ t\Phi(t^{-p}, t^{-1} \log t), & p^{-1} \in \N
    \end{cases}
    \quad \text{for} \;\; t \gg 1.
\end{equation}
This is exactly the desired expression \eqref{eq hyperbolic analytic function}. We complete the proof of Theorem~\ref{thm rigidity hyperbolic}.
\end{proof}

Finally, Corollary~\ref{cor hyperbolic} follows from the analytic representation \eqref{eq hyperbolic analytic function} after additional calculations of several leading coefficients.

\begin{proof}[Proof of Corollary~\ref{cor hyperbolic}]

Let $n$ be the largest nonnegative integer such that $np<1$. Expanding \eqref{eq analytic expansion of h} into power series, then all terms with exponential decay rates greater than $-1$ are of the form $c_j e^{-jps}$, where $1 \le j \le n$ and $c_j$ is a constant. We have
\begin{equation} \label{eq finite expansion of h preliminary}
    h(s)= \sum_{j=1}^n c_j e^{-jps}+ O(e^{-s}s) \quad \text{as} \;\; s \to +\infty.
\end{equation}
and the remainder can be strengthened to $O(e^{-s})$ in the case $p^{-1} \not\in \N$.

\begin{lem} \label{lem low order hyperbolic}
For $1 \le j \le n$, $c_j$ is uniquely determined by $v$ and is analytic in $v$.
\end{lem}

\begin{proof} 
Assume $n \ge 1$ and we proceed by induction on $j=1, \cdots, n$. Note this means $p<1$.

Set $F=\nabla U$, which is analytic on $\Y$. Then
\begin{equation} \label{eq equation of h hyperbolic}
    h''+h'= e^{-ps} F(v+h).
\end{equation}
By comparing the coefficients of $e^{-ps}$ in \eqref{eq equation of h hyperbolic}, we get 
\begin{equation}
    p(p-1) c_1= F(v).
\end{equation}
Since $p(p-1) \neq 0$, we see that $c_1= \frac{F(v)}{p(p-1)}$ is uniquely determined by $v$ and is analytic in $v$. 

Now suppose the conclusion holds for $c_1, \cdots, c_{j-1}$, where $2 \le j \le n$. Since $-jp>-1$, the coefficient of $e^{-jps}$ in $e^{-ps}F(v+h)$ is equal to that in
\begin{equation}
    e^{-ps} F \Big( v+\sum_{l=1}^{j-1} e^{-lps} c_l \Big).   
\end{equation}
By Taylor's theorem and the induction hypothesis, this coefficient is a constant $A_j$ depending only on $v$ and $j$. Comparing the coefficients of $e^{-jps}$ in \eqref{eq equation of h hyperbolic}, we obtain
\begin{equation}
    jp(jp-1)c_j= A_j.
\end{equation}
Again, since $jp(jp-1) \neq 0$, we see that $c_j= \frac{A_j}{jp(jp-1)}$ is uniquely determined by $v$ and $j$. This completes the induction.
\end{proof}

By \eqref{eq finite expansion of h preliminary} and Lemma~\ref{lem low order hyperbolic}, we also deduce
\begin{equation}
    e^{-ps} F(v+h)= \sum_{j=1}^{n+1} A_j(v) e^{-jps}+ O(e^{-(1+p)s}s) \quad \text{as} \;\; s \to +\infty,
\end{equation}
where $A_1, \cdots, A_{n+1}$ are analytic in $v$, and the remainder can be strengthened to $O(e^{-(1+p)s})$ in the case $p^{-1} \not\in \N$. Also note that all solutions of $y''+y'=0$ are generated by $1$ and $e^{-s}$, and recall $h(s) \to 0$ as $s \to +\infty$. Therefore, if $p^{-1} \not\in \N$, then 
\begin{equation}
    h(s)= \sum_{j=1}^n c_j e^{-jps}+ be^{-s}+ c_{n+1} e^{-(n+1)ps}+ O(e^{-(1+p)s}) \quad \text{as} \;\; s \to +\infty.
\end{equation}
If $p^{-1} \in \N$, then $(n+1)p=1$ and 
\begin{equation}
    h(s)= \sum_{j=1}^n c_j e^{-jps}+ c_{n+1} e^{-s} s+ be^{-s}+ O(e^{-(1+p)s}s) \quad \text{as} \;\; s \to +\infty.    
\end{equation}
By comparing the coefficients of $e^{-(n+1)ps}$ in \eqref{eq equation of h hyperbolic}, we have $(n+1)p ((n+1)p-1) c_{n+1}= A_{n+1}$ when $p^{-1} \not\in \N$ and $-c_{n+1}= A_{n+1}$ when $p^{-1} \in \N$. In both cases, we deduce that $c_{n+1}$ is uniquely determined by $v$ and is analytic in $v$.

Combining the above two expansions and $x(t)= vt+ th(t)$, we get
\begin{equation}
    x(t)= vt+ \sum_{j=1}^n c_j t^{1-jp}+ b+ c_{n+1} t^{1-(n+1)p}+ O(t^{-p}) \;\; \text{as} \;\; t \to +\infty \;\; (p^{-1} \not\in \N)
\end{equation} 
and
\begin{equation}
    x(t)= vt+ \sum_{j=1}^n c_j t^{1-jp}+ c_{n+1} \log t+ b+ O(t^{-p} \log t) \;\; \text{as} \;\; t \to +\infty \;\; (p^{-1} \in \N).
\end{equation} 
These are the desired expansions in Corollary~\ref{cor hyperbolic}. We have also shown that $c_1, \cdots, c_{n+1}$ depend analytically only on $v$, and hence are analytic functions from $\Y$ to $\X$.

The limits \eqref{eq limit defn of v} and \eqref{eq limit defn of b} are immediate. This defines the parametrization map $x \mapsto (v,b)$. Then \cite[Theorem~1]{Wu} asserts that the map is surjective. By the analyticity of $\Phi$, the expansion of $x$ yields 
\begin{equation}
    \dot{x}(t) = v+ \sum_{j=1}^{n+1} (1-jp)c_j t^{-jp}+ O(t^{-1-p}) \quad \text{as} \;\; t \to +\infty
\end{equation}
for $p^{-1} \not \in \N$, and
\begin{equation} 
    \dot{x}(t) = v+ \sum_{j=1}^n (1-jp)c_j t^{-jp}+ c_{n+1} t^{-1}+  O(t^{-1-p} \log t) \quad \text{as} \;\; t \to +\infty
\end{equation}
for $p^{-1} \in \N$. Then we can apply the uniqueness assertion in \cite[Theorem~1]{Wu} to obtain that the map is also injective. We thus finish the proof of Corollary~\ref{cor hyperbolic}. 
\end{proof}

\section{Variational results} \label{sec variational}

\subsection{Existence of minimizing profiled expansive solutions} \label{sec existence mixed} \

In this section, we give the proof of Theorem~\ref{thm existence mixed}. Assume $d \ge 2$ and $p \in (0,2)$. Let $(v,a)$ be an admissible pair and $x_0 \in \X$.

Let $r^0(t)= vt+ at^\frac{2}{p+2}+ f(t)$, where $f \in S_{1-p}$, be the approximate solution obtained in Lemma~\ref{lem approximate solution}. Choose $T_0>0$ large enough so that $r(t):=r^0(t+T_0) \in \Y$ for $t \ge 0$. Since $1-p<\frac{2}{p+2}$ for $p>0$, we have
\begin{equation}
    r(t) = vt+ at^{\frac{2}{p+2}}+ o(t^{\frac{2}{p+2}}) \quad \text{as} \;\; t \to +\infty
\end{equation}
Then by Lemma~\ref{lem lower bound on trajectory}, there exists $\kappa \in L^\infty \big( [0,+\infty); \R_{\ge 0} \big)$ satisfying \eqref{eq taylor lower bound on trajectory} and \eqref{eq limsup of kappa}. Since $\frac{p}{(p+2)^2}< \frac{1}{8}$ for $p \neq 2$, by increasing $T_0$ if necessary, we further have
\begin{equation} \label{eq upper bound of kappa}
    \sup_{t \ge 0} t^2 \kappa(t)< \frac{1}{8}.    
\end{equation}

For $0\le t_1< t_2< +\infty$ and $z \in H^1([t_1,t_2]; \X)$, we define the \emph{renormalized Lagrangian action functional}
\begin{equation}
    \A_{x_0}(z; [t_1,t_2]):= \int_{t_1}^{t_2} \Big[ \frac{1}{2} |\dot{z}|_\lambda^2+ U(r- r(0)+ x_0+ z)- U(r)- (\ddot{r}, z)_\lambda \Big] \dt.
\end{equation}
Set $x=r-r(0)+ x_0+z$. Then $\A_{x_0}$ and $\L$ are related by
\begin{equation} \label{eq relation of A and L}
    \A_{x_0}(z; [t_1,t_2])= \L(x; [t_1,t_2])- \L(r; [t_1,t_2])- {(\dot{r}, z)}_\lambda \big|_{t=t_1}^{t=t_2}.
\end{equation}
This follows from $|\dot{z}|_\lambda^2= |\dot{x}|_\lambda^2- |\dot{r}|_\lambda^2- 2(\dot{r}, \dot{z})_\lambda$ and integration by parts. 
\begin{rmk}
Note that $\L(r; [t_1,t_2])$ and ${(\dot{r}, z)}_\lambda \big|_{t=t_1}^{t=t_2}$ are always finite, but 
\begin{equation}
    \A_{x_0}(z; [t_1,t_2]), \; \L(x; [t_1,t_2]) \in \R \cup \{+\infty\}.
\end{equation}
Thus \eqref{eq relation of A and L} should be understood in the following sense: if either of these two quantities is $+\infty$, then so is the other; if either of them is finite, then the other is also finite and \eqref{eq relation of A and L} holds.    
\end{rmk}

We next extend the renormalized action to the half line $[0,+\infty)$. To this end, define
\begin{equation}
    \D = \Big\{ z \in H_{\text{loc}}^1 \big( [0,+\infty); \X \big) \ \big| \ z(0)=0, \ \dot{z} \in L^2 \big( [0,+\infty); \X \big) \Big\},
\end{equation}
equipped with the norm
\begin{equation}
    \|z\|_{\D}^2:= \int_0^\infty |\dot{z}(t)|_\lambda^2 \dt.
\end{equation}
The space $(\D, \|\cdot\|_{\D})$ is a Hilbert space. By the Cauchy--Schwarz inequality, we have
\begin{equation} \label{eq weighted L^infty}
    |z(t)|_\lambda \le \sqrt{t} \|z\|_{\D}, \quad \forall z \in \D, \ t \ge 0.
\end{equation}
Moreover, Hardy's inequality implies
\begin{equation} \label{eq hardy}
    \int_0^\infty t^{-2} |z(t)|_\lambda^2 \dt \le 4\|z\|_{\D}^2, \quad \forall z \in \D.
\end{equation} 

Then for $z \in \D$, we define the \emph{half-line renormalized Lagrangian action functional}
\begin{equation}
    \A_{x_0}(z):= \int_0^\infty \Big[ \frac{1}{2} |\dot{z}|_\lambda^2+ U(r- r(0)+ x_0+ z)- U(r)- (\ddot{r}, z)_\lambda \Big] \dt.
\end{equation}

\begin{lem}
The functional $\A_{x_0}$ satisfies the following properties.

\begin{enumerate}
    \item Coercivity. There exist $c, C>0$ independent of $z$ such that 
        \begin{equation}
            \A_{x_0}(z) \ge c\|z\|_{\D}^2 - C\|z\|_{\D}-C, \quad \forall z \in \D.
        \end{equation}
    \item Weak lower semicontinuity. If $z_n \rightharpoonup z$ in $\D$, then 
    \begin{equation}
        \liminf_{n \to +\infty} \A_{x_0}(z_n) \ge \A_{x_0}(z).
    \end{equation}
    \item Properness. There exists $z \in \D$ such that $\A_{x_0}(z)< +\infty$.
\end{enumerate}
\end{lem}

\begin{proof}
Set $y_0= x_0- r(0)$.

(1) Let $c, C>0$ denote constants independent of $t$ and $z$. We decompose
\begin{equation} 
    \A_{x_0}(z) = F(z)+ Q(z)- L(z)+ K,
\end{equation}
where
\begin{align}
    F(z) &= \int_0^\infty \Big[ U(r+ y_0+ z)- U(r)- (\nabla U(r), z+y_0)_\lambda+ \kappa |z+y_0|_\lambda^2 \Big] \dt, \\
    Q(z) &= \frac{1}{2} \|z\|_{\D}^2- \int_0^\infty \kappa|z|_\lambda^2 \dt, \\
    L(z) &= \int_0^\infty \big( \ddot{r}- \nabla U(r)+ 2\kappa y_0, z \big)_\lambda \dt, \\
    K &= \int_0^\infty \Big[ (\nabla U(r), y_0)_\lambda - \kappa |y_0|_\lambda^2 \Big] \dt. 
\end{align}
Here $\kappa$ is chosen at the beginning of this section. We have \eqref{eq taylor lower bound on trajectory}, \eqref{eq upper bound of kappa} and $\kappa \in L^\infty([0,+\infty))$.

By \eqref{eq taylor lower bound on trajectory}, we deduce $F(z) \ge 0$. Using \eqref{eq upper bound of kappa} and Hardy's inequality \eqref{eq hardy}, we get 
\begin{equation} \label{eq coercivity term1} \begin{aligned}
    Q(z) &\ge \frac{1}{2} \|z\|_{\D}^2-  \sup_{t \ge 0} t^2 \kappa(t) \int_0^\infty t^{-2} |z|_\lambda^2 \dt \\
    &\ge \Big( \frac{1}{2}-4\sup_{t \ge 0} t^2 \kappa(t) \Big) \|z\|_{\D}^2 \ge c\|z\|_{\D}^2.
\end{aligned} \end{equation}
By \eqref{eq approximate solution estimate}, we have $|\ddot{r}- \nabla U(r)| \le C(t+1)^{-2}$. We also have $\kappa(t) \le C(t+1)^{-2}$ by $\kappa \in L^\infty([0,+\infty))$ and \eqref{eq upper bound of kappa}. Thus by \eqref{eq weighted L^infty}, we obtain
\begin{equation} \label{eq coercivity term2} \begin{aligned}
    |L(z)| &\le C\int_0^\infty (t+1)^{-2} |z|_\lambda \dt \\
    &\le C \int_0^\infty \sqrt{t} (t+1)^{-2} \dt \ \|z\|_{\D} \le C \|z\|_{\D}.
\end{aligned} \end{equation}
Since $|r_i-r_j| \ge c(t+1)^{\frac{2}{p+2}}$ for $i \neq j$, one has $|\nabla U(r)| \le C(t+1)^{-\frac{2(p+1)}{p+2}}$. Thus
\begin{equation} \label{eq coercivity term3} \begin{aligned}
    |K| \le C\int_0^\infty \Big[ (t+1)^{-\frac{2(p+1)}{p+2}} + (t+1)^{-2} \Big] \dt \le C.
\end{aligned} \end{equation}
Gathering \eqref{eq coercivity term1}, \eqref{eq coercivity term2}, \eqref{eq coercivity term3} and $F(z) \ge 0$ together, we deduce
\begin{equation}
    \A_{x_0}(z) \ge c \|z\|_{\D}^2- C\|z\|_{\D}- C.
\end{equation}
This proves coercivity. It also implies  $\A_{x_0}: \D \to \R \cup \{+\infty\}$.

(2) Assume $z_n \rightharpoonup z$ in $\D$. Note that \eqref{eq coercivity term1} implies $Q(z)$ is a bounded positive quadratic form on $\D$, and \eqref{eq coercivity term2} implies $L(z)$ is a bounded linear functional on $\D$. Thus
\begin{equation}
    \liminf_{n \to +\infty} Q(z_n) \ge Q(z), \qquad \lim_{n \to +\infty} L(z_n)= L(z).
\end{equation}
Moreover, $z_n(t) \to z(t)$ for all $t \ge 0$. Since the integrand in $F$ is non-negative, by lower semicontinuity of $U$ and Fatou's lemma,
\begin{equation}
    \liminf_{n \to +\infty} F(z_n) \ge F(z).
\end{equation}
Then we conclude $\liminf_{n \to +\infty} \A_{x_0}(z_n) \ge \A_{x_0}(z)$.

(3) Since $p \in (0,2)$, there exists $\beta \in (1/2, 1/p)$. Take $y \in \Y$, $\epsilon>0$ and set
\begin{equation}
    x(t)= x_0+ t^\beta y \quad \text{for}  \;\; 0 \le t \le \epsilon.
\end{equation}
If $\epsilon$ is small enough, then $x(t) \in \Y$ for $0< t \le \epsilon$. In particular, $x(\epsilon) \in \Y$.

Next, choose $T$ large enough so that $r(t)+ y_0 \in \Y$ for all $t \ge T$. Since $d \ge 2$, $\Y$ is path connected, so we can extend $x$ on $[\epsilon, T]$ as a continuous and piecewise smooth curve in $\Y$, with $x(T)= r(T)+ y_0$. Finally, for $t>T$, we take $x(t)= r(t)+ y_0$.

Now we have defined $x \in C \big( [0, +\infty); \X \big)$. Let $z(t)= x(t)- r(t)- y_0$. We claim that this $z$ satisfies $z \in \D$ and $\A_{x_0}(z)< +\infty$.

In fact, our construction guarantees $z(0)=0$ and $z(t)=0$ for $t \ge T$, so $z$ is compactly supported. Since $z$ is continuous, $z \in L^2 \big( [0,+\infty); \X \big)$. By the definition of $x$ on $[0, \epsilon]$, we have
\begin{equation}
    |\dot{x}|^2= O(t^{2\beta-2}), \quad U(x)= O(t^{-p\beta}) \quad \text{as} \;\; t \to 0^+.
\end{equation}
Notice $2\beta-2> -1$ and $-p\beta> -1$ by our choice of $\beta$. It then follows that $\dot{z} \in L^2 \big( [0,+\infty); \X \big)$ and $\int_0^T U(x(t)) \dt< +\infty$. Thus $z \in \D$. For $t>T$, we have
\begin{equation}
    U(x(t))- U(r(t))= \big( \nabla U(r(t)+ \theta y_0), y_0 \big)_\lambda \quad \text{for some} \;\;\theta \in [0,1].
\end{equation}
We have $|\nabla U(r(t)+ \theta y_0)| \le Ct^{-\frac{2(p+1)}{p+2}}$ because $|r_i-r_j| \ge ct^{\frac{2}{p+2}}$, and thus
\begin{equation}
    \int_T^\infty \big| U(x(t))- U(r(t)) \big| \dt< +\infty.
\end{equation}
We conclude that $\A_{x_0}(z)< +\infty$.
\end{proof}

Since $\D$ is a Hilbert space, with the above lemma, the direct method in the calculus of variations (see \cite[Chapter~I, Theorem~1.2]{Struwe}) implies that $\A_{x_0}$ admits a minimizer $z_* \in \D$. By the properness of $\A_{x_0}$, we have 
\begin{equation} \label{eq finite action}
    \A_{x_0}(z_*)< +\infty \;\; \text{and} \;\; \A_{x_0}(z_*; [t_1,t_2])< +\infty, \;\; \forall 0 \le t_1< t_2< +\infty.
\end{equation}

Now we can give the proof of Theorem~\ref{thm existence mixed}.

\begin{proof}[Proof of Theorem~\ref{thm existence mixed}]  
Let $x_*(t)= r(t)+ z_*(t)+ y_0$, where $y_0= x_0- r(0)$ as before.

We have $x_*(0)=x_0$ because $z_*(0)=0$. By \eqref{eq weighted L^infty}, we have $z_*(t)= O(t^{\frac{1}{2}})$. Also recall 
\begin{equation}
    r(t)= v(t+T_0)+ a(t+T_0)^{\frac{2}{p+2}}+ f(t+T_0)
\end{equation} 
with $f \in S_{1-p}$. Note $p \in (0,2)$ implies $1-p< \frac{2}{p+2}$ and $\frac{1}{2}< \frac{2}{p+2}$. We thus have 
\begin{equation} \label{eq minimizing solution expansion}
    x_*(t)= vt+ at^{\frac{2}{p+2}}+ o(t^{\frac{2}{p+2}}) \quad \text{as} \;\; t \to +\infty.
\end{equation}

Next, we show that $x_*$ is a solution of \eqref{eq N-body}. By Remark~\ref{rmk Marchal}, it suffices to show that $x_*$ is a fixed-time minimizer on any interval $[t_1,t_2]$, where $0<t_1<t_2<+\infty$. 

Using \eqref{eq finite action} and that $r$ is collision-free, we have $\L(x_*; [t_1,t_2])< +\infty$. Let $x \in H^1([t_1,t_2]; \X)$ be such that $x(t_1)= x_*(t_1)$ and $x(t_2)= x_*(t_2)$. We may assume $\L(x; [t_1,t_2])< +\infty$, since otherwise $\L(x; [t_1,t_2]) \ge \L(x_*; [t_1,t_2])$ trivially. We define
\begin{equation}
    z(t)= \begin{cases}
        z_*(t), & t \not\in[t_1,t_2], \\
        x(t)- r(t)- y_0, & t \in [t_1,t_2].
    \end{cases}
\end{equation}
Then $z \in \D$. Moreover, since $z=z_*$ outside $(t_1,t_2)$ and \eqref{eq finite action}, we have
\begin{equation}
    \A_{x_0}(z)- \A_{x_0}(z_*)= \A_{x_0}(z; [t_1,t_2])- \A_{x_0}(z_*; [t_1,t_2]).
\end{equation}
By \eqref{eq relation of A and L} and $z(t_1)=z_*(t_1)$, $z(t_2)= z_*(t_2)$, we deduce
\begin{equation}
    \A_{x_0}(z)- \A_{x_0}(z_*)= \L(x; [t_1,t_2])- \L(x_*; [t_1,t_2]).
\end{equation}
This implies $\L(x; [t_1,t_2]) \ge \L(x_*; [t_1,t_2])$ by the minimality of $z_*$. Thus $x_*$ solves \eqref{eq N-body}.

Finally, we prove the free-time minimality of $x_*$ at energy $E= \f12 |v|_\lambda^2$. 

Fix $0\le t_1< t_2< +\infty$. We have $\L(x_*; [t_1,t_2])< +\infty$ by \eqref{eq finite action} and since $r(t)$ is collision-free for all $t \ge 0$. Let $x \in H^1([0, T]; \X)$ with $T>0$, $x(0)= x_*(t_1)$ and $x(T)= x_*(t_2)$. We may assume $\L(x; [0, T])< +\infty$, since otherwise \eqref{eq free time minimizer} holds trivially. We construct $\hat{x}$ by connecting three curves $x_* |_{[0,t_1]}$, $x|_{[0,T]}$ and $x_* |_{[t_2,+\infty)}$. More precisely, we set $c=t_2-t_1-T$ and define
\begin{equation}
    \hat{x}(t)= \begin{cases}
        x_*(t), &0 \le t \le t_1, \\
        x(t-t_1), &t_1 \le t \le t_1+T, \\
        x_*(t+c), &t \ge t_1+T. 
    \end{cases}
\end{equation}
Let $z= \hat{x}- r- y_0$. Note $z(t)= z_*(t)+ x_*(t+c)-x_*(t)$ for $t$ large, and 
\begin{equation} \begin{aligned}
    \dot{x}_*(t+c)- \dot{x}_*(t)= \int_t^{t+c} \ddot{x}_*(s) \d s&= \int_t^{t+c} \nabla U(x_*(s)) \d s \\
    &= O(t^{-\frac{2(p+1)}{p+2}}) \quad \text{as} \;\; t \to +\infty.
\end{aligned} \end{equation}
We thus have $\dot{z} \in L^2([0,+\infty); \X)$. We deduce $z \in \D$. 

By \eqref{eq relation of A and L} and \eqref{eq finite action}, for any $R>0$, we have
\begin{equation} \label{eq identity for free time minimizer} \begin{aligned}
    &\quad \ \A_{x_0}(z; [0,R])- \A_{x_0}(z_*; [0,R]) \\
    &= \L(\hat{x}; [0,R])- \L(x_*; [0,R])- {(\dot{r}, \hat{x}-x_*)}_\lambda \big|_{t=0}^{t=R}.
\end{aligned} \end{equation}
For $R>\max \{t_1+T, t_2\}$, we have
\begin{equation} \begin{aligned}
    \L(\hat{x}; [0,R])- \L(x_*; [0,R]) &= \L(x; [0,T])- \L(x_*; [t_1,t_2]) \\
    &\quad \ + \int_R^{R+c} \Big[ \frac{1}{2} |\dot{x}_*(t)|_\lambda^2+ U(x_*(t)) \Big] \dt.
\end{aligned} \end{equation}
Since $\dot{x}_*(t) \to v$ and $U(x_*(t)) \to 0$ as $t \to +\infty$, we deduce
\begin{equation}
    \lim_{R \to +\infty} \int_R^{R+c} \Big[ \frac{1}{2} |\dot{x}_*(t)|_\lambda^2+ U(x_*(t)) \Big] \dt= \f12 |v|_\lambda^2 c= Ec.
\end{equation}
Besides, using $\hat{x}(0)= x_*(0)$, for $R>t_1+T$, we have
\begin{equation} \begin{aligned}
    {(\dot{r}, \hat{x}-x_*)}_\lambda \big|_{t=0}^{t=R} &= \big( \dot{r}(R), \hat{x}(R)- x_*(R) \big)_\lambda \\
    &= \big( \dot{r}(R), x_*(R+c)- x_*(R) \big)_\lambda.
\end{aligned} \end{equation}
We have $\dot{r}(R) \to v$ as $R \to +\infty$ by the expansion of $r$ and $f \in S_{1-p}$, while $\dot{x}_* \to v$ implies
\begin{equation}
    x_*(R+c)- x_*(R)= \int_R^{R+c} \dot{x}_*(t) \dt \to cv \quad \text{as} \;\; R \to +\infty.
\end{equation}
Thus
\begin{equation}
    \lim_{R \to +\infty} \big( \dot{r}(R), x_*(R+c)- x_*(R) \big)_\lambda= (v, cv)_\lambda= 2E c.
\end{equation}
Combining these limits, by taking $R \to +\infty$ in \eqref{eq identity for free time minimizer}, we get
\begin{align}
    \A_{x_0}(z)- \A_{x_0}(z_*) &= \L(x; [0,T])- \L(x_*; [t_1,t_2])+ Ec- 2Ec \\
    &= \L(x; [0,T])+ ET- \L(x_*; [t_1,t_2])-E(t_2-t_1).
\end{align}
We deduce \eqref{eq free time minimizer} by the minimality of $z_*$.

The proof of Theorem~\ref{thm existence mixed} is finished.
\end{proof}

\subsection{Existence of hyperbolic solutions in the strong-force regime} \label{sec existence hyperbolic} \

In this section, we prove Theorem~\ref{thm existence hyperbolic}. Assume $d \in \N_{\ge 2}$ and $p \ge 2$. Let $v \in \Y$ and $x_0 \in \X$.

The main difficulty is that the functional $\A_{x_0}$ defined in Section~\ref{sec existence mixed} fails to be proper when $x_0 \in \Delta$. In fact, the Lagrangian action functional $\L$ must blow up near $t=0$. To see this, suppose $(x_0)_i= (x_0)_j$ for some $i \neq j$. Then for any $T>0$ and $x \in H^1([0,T])$ with $x(0)= x_0$, we have 
\begin{equation}
    x_i(t)- x_j(t)= \int_0^t \big( \dot{x}_i(\tau)- \dot{x}_j(\tau) \big) \d\tau.
\end{equation}
By the Cauchy--Schwarz inequality,
\begin{equation}
    |x_i(t)- x_j(t)| \le \int_0^t \big| \dot{x}_i(\tau)- \dot{x}_j(\tau) \big| \d\tau \le C\|x\|_{H^1} \sqrt{t}.
\end{equation}
Using $p \ge 2$, we see
\begin{equation}
    \L(x; [0,T]) \ge \int_0^T U(x(t)) \dt \ge \int_0^T (C\|x\|_{H^1})^{-p} t^{-\frac{p}{2}} \dt= +\infty.
\end{equation}
The above calculation shows that the same proof as in Section~\ref{sec existence mixed} cannot be directly applied when $x_0 \in \Delta$. On the other hand, when $x_0 \in \Y$, that argument can be adapted. 

To overcome this difficulty, we construct a renormalized action functional that not only captures the asymptotic behavior at infinity time, but also fits the initial collision dynamics. The following proof applies to arbitrary $x_0 \in \X$, regardless of whether $x_0$ is collisional.

We make the ansatz $x(t)= x_0+ t^{\frac{2}{p+2}} y(t)$ for small $t$. Write $U= U_0+ U_1$, where
\begin{equation}
    U_0= \sum_{I \in M_{x_0}} U_I, \qquad U_1= \frac{1}{2} \sum_{I \neq J \in M_{x_0}} U_{IJ}.
\end{equation}
Then we can compute the Lagrangian of $x$
\begin{equation} \label{eq one step calculation} \begin{aligned}
    \f12 |\dot{x}|_\lambda^2+ U(x) &= \f12 t^{\frac{4}{p+2}} |\dot{y}|_\lambda^2+ t^{-\frac{2p}{p+2}} W(y) \\
    &\quad \ + U_1(x_0+ t^{\frac{2}{p+2}} y) + \frac{1}{p+2} \f\d\dt \big( t^{\frac{2-p}{p+2}} |y|_\lambda^2 \big),
\end{aligned} \end{equation}
where 
\begin{equation}
     W(y)= U_0(y)+ \frac{p}{(p+2)^2} |y|_\lambda^2.
\end{equation}
For $a \in \X$ such that $(x_0,a)$ is an admissible pair, we have \eqref{eq Taylor around a_I} and $\nabla U_I(a_I)= -\frac{2p}{(p+2)^2} a_I$ for every $I \in M_{x_0}$. Summing over $I$, we get 
\begin{equation}
    W(a)= \min_{y \in \X} W(y) \triangleq m. 
\end{equation}

The above calculation inspires the following definitions. Let
\begin{equation}
    \H= \Big\{y \in H_{\text{loc}}^1 \big( (0,1]; \X \big) \ \big| \ t^{\frac{2}{p+2}} \dot{y} \in L^2 \big( (0,1]; \X \big) \Big\}
\end{equation}
equipped with the norm
\begin{equation}
    \|y\|_{\H}^2:= \int_0^1 t^\frac{4}{p+2} |\dot{y}(t)|_\lambda^2 \dt+ |y(1)|_\lambda^2,
\end{equation}
and
\begin{equation}
    \D_1 = \Big\{ z \in H_{\text{loc}}^1 \big( [1,+\infty); \X \big) \ \big| \ z(1)=0, \ \dot{z} \in L^2 \big( [1,+\infty); \X \big) \Big\}
\end{equation}
equipped with the norm
\begin{equation}
    \|z\|_{\D_1}^2:= \int_1^\infty |\dot{z}(t)|_\lambda^2 \dt.
\end{equation}
Note that $\D_1$ is simply a time-translated version of $\D$ in Section~\ref{sec existence mixed}. The spaces $(\H, \| \cdot \|_{\H})$ and $(\D_1, \| \cdot \|_{\D_1})$ are both Hilbert spaces. By the Cauchy--Schwarz inequality, for any $y \in \H$ and $ t \in (0,1]$, we have
\begin{equation} \label{eq pointwise estimate in H}
    |y(t)|_\lambda \le \begin{cases}
        |y(1)|_\lambda+ C(p) \left( \int_0^1 t^\frac{4}{p+2} |\dot{y}(t)|_\lambda^2 \dt \right)^{1/2}, & p>2, \\
        |y(1)|_\lambda+ \sqrt{|\log t|} \left( \int_0^1 t^\frac{4}{p+2} |\dot{y}(t)|_\lambda^2 \dt \right)^{1/2}, & p=2.
    \end{cases}
\end{equation}

For $(y,z) \in \H \times \D_1$, we define
\begin{align}
    \A_{x_0}(y,z) &= \int_0^1 \Big[ \f12 t^{\frac{4}{p+2}} |\dot{y}|_{\lambda}^2+ t^{-\frac{2p}{p+2}} (W(y)-m)+ U_1(x_0+ t^{\frac{2}{p+2}} y) \Big] \dt \\
    &\quad \ + \frac{1}{p+2} |y(1)|_\lambda^2 \\
    &\quad \ + \int_1^\infty \Big[ \frac{1}{2} |\dot{z}|_\lambda^2+ U(vt+ b_y+ z)- U(vt) \Big] \dt \\
    &\quad \ - (v, b_y)_\lambda,
\end{align}
where $b_y:= x_0+ y(1)- v$. Here, the rescaled variable $y$ describes the behavior near the initial collision, while $z$ controls the hyperbolic behavior near infinity. They are coupled through the matching constant $b_y$.

\begin{lem}
$\A_{x_0}$ is coercive, weakly lower semicontinuous and proper on $\H \times \D_1$.
\end{lem}

\begin{proof}
Coercivity of $\A_{x_0}$ is clear. In fact, by $U \ge U_1 \ge 0$ and $W \ge m$, we have
\begin{equation} \begin{aligned}
    \A_{x_0}(y,z) &\ge \f12 \int_0^1 t^{\frac{4}{p+2}} |\dot{y}|_{\lambda}^2 \dt + \frac{1}{p+2} |y(1)|_\lambda^2 \\
    &\quad \ + \f12 \int_1^\infty |\dot{z}|_\lambda^2- \int_1^\infty U(vt) \dt- (v, b_y)_\lambda.
\end{aligned} \end{equation}
Using $|b_y| \le |y(1)|+ C$ and $U(vt)= t^{-p} U(v)$, which is integrable on $[1,+\infty)$, we get 
\begin{equation}
    \A_{x_0}(y,z) \ge c\|y\|_\H^2+ c\|z\|_{\D_1}^2- C\|y\|_{\H}- C.
\end{equation}

The kinetic terms are weakly lower semicontinuous since they are positive definite quadratic forms. For potential terms, by \eqref{eq pointwise estimate in H} and a time-translated version of \eqref{eq weighted L^infty}, we see that $(y_n,z_n) \rightharpoonup (y,z)$ in $\H \times \D_1$ implies $y_n(t) \to y(t)$ for each $t \in (0,1]$ and $z_n(t) \to z(t)$ for each $t \in [1, +\infty)$. Then by lower semicontinuity of $U$ and Fatou's lemma, we obtain weak lower semicontinuity.

Finally, we show that $\A_{x_0}$ is proper. Take $a \in \X$ such that $(x_0, a)$ is an admissible pair. Then choose $\epsilon \ll 1$ such that $x(t):= x_0+ at^{\frac{2}{p+2}} \in \Y$ for all $t \in (0, \epsilon]$. By the path connectivity of $\Y$ when $d \ge 2$, we can extend $x$ on $[\epsilon, 1]$ such that $x \in C^2((0,1]; \Y)$ and $x(1)=v$. Take
\begin{equation}
    y(t)= t^{-\frac{2}{p+2}} (x(t)- x_0), \quad \forall 0< t \le 1 \quad \text{and} \quad z \equiv 0.
\end{equation}
Then we have $y \in \H$ and $z \in \D_1$. We claim that $\A_{x_0}(y,z)< +\infty$ holds. Since $y(1)= v-x_0$ and thus $b_y=0$, it suffices to show 
\begin{equation}
    \int_0^1 \Big[ t^{-\frac{2p}{p+2}} (W(y)-m)+ U_1(x_0+ t^{\frac{2}{p+2}} y) \Big] \dt< +\infty.
\end{equation}
Note $y(t)= a$ for $t \in (0,\epsilon]$, so $W(y)=m$ on $(0, \epsilon]$. We also have $U_1(x_0+ a t^{\frac{2}{p+2}}) \le C$ for $t \in (0, \epsilon]$ because $U_1$ only contains inter-cluster terms. Thus, the integral on $(0, \epsilon]$ is finite. On $[\epsilon, 1]$, we have $x \in \Y$ and the integrand is continuous, so the integral is also finite.
\end{proof}

\begin{proof}[Proof of Theorem~\ref{thm existence hyperbolic}]

We can now apply the direct method \cite[Chapter~I, Theorem~1.2]{Struwe} to conclude that $\A_{x_0}$ admits a minimizer $(y_*,z_*) \in \H \times \D_1$, and $\A_{x_0}(y_*,z_*)< +\infty$. Define
\begin{equation}
    x_*(t)= \begin{cases}
        x_0, & t=0, \\
        x_0+ t^{\frac{2}{p+2}} y_*(t), & 0 < t \le 1, \\
        vt+ b_{y_*}+ z_*(t), & t \ge 1.
    \end{cases}
\end{equation}
Since $b_{y_*}= x_0+ y_*(1)-v$, the above definition agrees at $t=1$. Moreover, using \eqref{eq pointwise estimate in H}, we see that $x_*(t)$ is continuous at $t=0$. Therefore, 
\begin{equation}
    x_* \in H_{\text{loc}}^1((0,+\infty); \X)\cap C^0([0,+\infty); \X) \quad \text{and} \quad x_*(0)=x_0.
\end{equation}
By \eqref{eq weighted L^infty}, we have $z_*(t)= O(\sqrt{t})$, and thus $x(t)=vt+ o(t)$ as $t \to +\infty$. 

It remains to show that $x_*$ solves \eqref{eq N-body} for $t>0$. By Remark~\ref{rmk Marchal}, it suffices to show that $x_*$ is a fixed-time minimizer on any interval $[t_1, t_2]$, where $0<t_1<t_2<+\infty$. We may only consider the case where $t_1<1<t_2$ by extending the interval.

By the definition of $\A_{x_0}$ and $x_*$, we see that $\A_{x_0}(y_*,z_*)< +\infty$ implies $\L(x_*; [t_1, t_2])<+\infty$. Let $x \in H^1([t_1,t_2]; \X)$ with $x(t_1)= x_*(t_1)$ and $x(t_2)=x_*(t_2)$. We may assume $\L(x;[t_1,t_2])<+\infty$. Take 
\begin{equation}
    y(t)= \begin{cases}
        y_*(t), & 0< t \le t_1, \\
        t^{-\frac{2}{p+2}} (x(t)- x_0), & t_1 \le t \le 1
    \end{cases}
\end{equation}
and
\begin{equation}
    z(t)= \begin{cases}
        x(t)- vt- b_y, & 1 \le t \le t_2, \\
        z_*(t)+ b_{y_*}- b_y, & t_2 \le t < +\infty.
    \end{cases}
\end{equation}
The definitions of $y$ and $z$ agree at $t_1$ and $t_2$, respectively, by the endpoint conditions of $x$. Moreover, $z(1)=0$ by the definition of $b_y$. Using $x \in H^1([t_1, t_2]; \X)$, we get $y \in H^1([t_1,1]; \X)$ and $z \in H^1([1,t_2]; \X)$. Combining these with $y=y_*$ on $(0,t_1]$ and $\dot{z}= \dot{z}_*$ on $[t_2,+\infty)$, we deduce $(y,z) \in \H \times \D_1$. We shall use $(y,z)$ as a competitor of $(y_*,z_*)$.

By $y=y_*$ on $(0, t_1]$ and $vt+ b_y+z= vt+ b_{y_*}+ z_*$ on $[t_2, +\infty)$, we have
\begin{equation}
    \A_{x_0}(y,z)- \A_{x_0}(y_*,z_*)= \L_1(y)- \L_1(y_*)+ \L_2(y, z)- \L_2(y_*, z_*),
\end{equation}
where
\begin{equation} \begin{aligned}
    \L_1(y) &= \int_{t_1}^1 \Big[ \f12 t^{\frac{4}{p+2}} |\dot{y}|_{\lambda}^2+ t^{-\frac{2p}{p+2}} W(y)+ U_1(x_0+ t^{\frac{2}{p+2}} y) \Big] \dt \\
    &\quad \ + \frac{1}{p+2} |y(1)|_{\lambda}^2
\end{aligned} \end{equation}
and 
\begin{equation}
    \L_2(y,z)= \int_1^{t_2} \Big[ \f12 |\dot{z}|_\lambda^2+ U(vt+b_y+z)- U(vt) \Big] \dt- (v,b_y)_\lambda.
\end{equation}

Since $x(t)=x_0+ t^{\frac{2}{p+2}} y(t)$ for $t_1 \le t \le 1$, using \eqref{eq one step calculation}, we have
\begin{equation}
    \L_1(y)= \L(x; [t_1,1])+ \frac{1}{p+2} t_1^{\frac{2-p}{p+2}} |y(t_1)|_\lambda^2.
\end{equation}
Since $x(t)= vt+ b_y+ z(t)$ and $|\dot{z}|_\lambda^2= |\dot{x}|_\lambda^2- |v|_\lambda^2- 2(v, \dot{z})_\lambda$ for $1 \le t \le t_2$, and $z(1)=0$, we have
\begin{equation}
    \L_2(y,z)= \L(x; [1,t_2])- \L(vt; [1,t_2])- (v, x(t_2)- vt_2)_\lambda.
\end{equation}
Then we get
\begin{equation} \begin{aligned}
    \L_1(y)+ \L_2(y,z) &= \L(x; [t_1, t_2])+ \frac{1}{p+2} t_1^{\frac{2-p}{p+2}}|y(t_1)|_\lambda^2 \\
    &\quad \ - \L(vt; [1,t_2])- (v, x(t_2)- vt_2)_\lambda.
\end{aligned} \end{equation}
The same identity holds for $x_*$, $y_*$ and $z_*$. Then by $y(t_1)= y_*(t_1)$ and $x(t_2)= x_*(t_2)$, we have
\begin{equation}
    \A_{x_0}(y,z)- \A_{x_0}(y_*,z_*)= \L(x; [t_1,t_2])- \L(x_*; [t_1,t_2]).
\end{equation}
Using the minimality of $(y_*,z_*)$, we get $\L(x; [t_1,t_2]) \ge \L(x_*; [t_1,t_2])$. This proves the fixed-time minimality of $x_*$ on any $[t_1,t_2] \subset (0,+\infty)$.

Finally, the free-time minimality of $x_*$ on $[0,+\infty)$ at energy $E= \f12 |v|_\lambda^2$ follows from the same argument as in Section~\ref{sec existence mixed}. If $x_0 \in \Delta$, then any curve has infinite action near $t=0$, which is allowed by our definition. Otherwise, we insert the competitor into $x_*$ by translating the subsequent trajectory in time. Then we can construct a competitor of $(y_*, z_*)$ and use an analogous action comparison between $\A_{x_0}$ and $\L$ to conclude free-time minimality.

We complete the proof of Theorem~\ref{thm existence hyperbolic}.
\end{proof}

\subsection{Eventual free-time minimality} \

Theorem~\ref{thm eventual minimizer} is a consequence of the global lower bound Lemma~\ref{lem lower bound on trajectory} and the variational framework in Section~\ref{sec existence mixed}. Assume $d \in \N_{\ge 1}$ and $p>0,\ p \neq 2$. Let $(v,a)$ be an admissible pair. Assume $x \in C^2((T,+\infty); \Y)$ solves \eqref{eq N-body} and satisfies \eqref{eq profiled expansive solution}.

\begin{proof}[Proof of Theorem~\ref{thm eventual minimizer}] 
For $T_0>T$, let $r(t)= x(t+T_0)$. Then $r$ satisfies \eqref{eq profiled expansive solution} as well. We can apply Lemma~\ref{lem lower bound on trajectory} to obtain a function $\kappa \in L^\infty([0,+\infty); \R_{\ge 0})$  satisfying \eqref{eq taylor lower bound on trajectory} and \eqref{eq limsup of kappa}. Since $\frac{p}{(p+2)^2}< \f18$ for $p \neq 2$, by increasing $T_0$ if necessary, we may replace \eqref{eq limsup of kappa} by $\sup_{t \ge 0} t^2 \kappa(t)< \f18$.

Let $\D$ be as in Section~\ref{sec existence mixed} and for $z \in \D$, define
\begin{equation}
    \A(z):= \int_0^\infty \Big[ \frac{1}{2} |\dot{z}|_\lambda^2+ U(r+z)- U(r)- (\ddot{r}, z)_\lambda \Big] \dt.
\end{equation}
Since $r$ is an exact solution of \eqref{eq N-body}, we have $\ddot{r}= \nabla U(r)$. Then by \eqref{eq taylor lower bound on trajectory} and \eqref{eq hardy}, we get
\begin{equation}
    \A(z) \ge \int_0^\infty \Big[ \frac{1}{2} |\dot{z}|_\lambda^2- \kappa |z|_\lambda^2 \Big] \dt \ge \Big( \f12- 4\sup_{t \ge 0} t^2 \kappa(t) \Big) \|z\|_{\D}^2 \ge c \|z\|_{\D}^2.
\end{equation}
Note also $\A(0)=0$. We see that $z_*=0$ is the unique minimizer of $\A$. 

By the proof of Theorem~\ref{thm existence mixed}, we deduce that $r=r+z_*$ is a free-time minimizer at energy $|v|_\lambda^2/2$. This means $x|_{[T_0,+\infty)}$ is a free-time minimizer at energy $|v|_\lambda^2/2$.
\end{proof}

\end{document}